\documentclass[11pt]{amsart}

\usepackage{bm}
\usepackage{amsxtra}
\usepackage{amssymb}
\usepackage{mathtools}
\usepackage{color}
\usepackage{enumitem}
\usepackage{tikz}
\usepackage{tikz-cd}
\usetikzlibrary{arrows}
\usepackage[all]{xy}
\usepackage{xcolor}
\usepackage{soul}
\usepackage[normalem]{ulem}
\usepackage{amsmath}
\newcommand{\stkout}[1]{\ifmmode\text{\sout{\ensuremath{#1}}}\else\sout{#1}\fi}
\newcommand\redsout{\bgroup\markoverwith{\textcolor{red}{\rule[0.5ex]{2pt}{0.4pt}}}\ULon}

\usetikzlibrary{arrows.meta}
\theoremstyle{plain}
\newcommand{\cleqn}{\setcounter{equation}{0}}
\newcommand{\clth}{\setcounter{theorem}{0}}
\newcommand {\sectionnew}[1]{\section{#1}\cleqn\clth}
\newcommand{\nn}{\hfill\nonumber}
\newtheorem{theorem}{Theorem}[section]
\newtheorem{lemma}[theorem]{Lemma}
\newtheorem{definition-theorem}[theorem]{Definition-Theorem}
\newtheorem{proposition}[theorem]{Proposition}
\newtheorem{corollary}[theorem]{Corollary}
\newtheorem{definition}[theorem]{Definition}
\newtheorem{construction}[theorem]{Construction}
\newtheorem{example}[theorem]{Example}
\newtheorem{remark}[theorem]{Remark}
\newtheorem{conjecture}[theorem]{Conjecture}

\newcommand \bth[1] { \begin{theorem}\label{t#1} }
\newcommand \ble[1] { \begin{lemma}\label{l#1} }

\newcommand \bpr[1] { \begin{proposition}\label{p#1} }
\newcommand \bco[1] { \begin{corollary}\label{c#1} }
\newcommand \bcon[1] { \begin{construction}\label{co#1} }
\newcommand \bde[1] { \begin{definition}\label{d#1}\rm }
\newcommand \bex[1] { \begin{example}\label{e#1}\rm }
\newcommand \bre[1] { \begin{remark}\label{r#1}\rm }
\newcommand \bcj[1] { \begin{conjecture}\label{j#1}\rm }

\renewcommand {\eth} { \end{theorem} }
\newcommand {\ele} { \end{lemma} }

\newcommand {\epr} { \end{proposition} }
\newcommand {\eco} { \end{corollary} }
\newcommand {\econ} { \end{construction} }
\newcommand {\ede} { \end{definition} }
\newcommand {\eex} { \end{example} }
\newcommand {\ere} { \end{remark} }
\newcommand {\ecj} { \end{conjecture} }

\newcommand {\enota} { \end{notation} }
\newcommand \thref[1]{Theorem \ref{t#1}}
\newcommand \leref[1]{Lemma \ref{l#1}}
\newcommand \prref[1]{Proposition \ref{p#1}}
\newcommand \coref[1]{Corollary \ref{c#1}}

\newcommand \deref[1]{Definition \ref{d#1}}
\newcommand \exref[1]{Example \ref{e#1}}
\newcommand \reref[1]{Remark \ref{r#1}}
\newcommand \lb[1]{\label{#1}}
\def \Abe {{\mathcal A}_\varepsilon}
\def \Abbe {{\mathcal A}_\varepsilon^{1/2}}
\def \AA {{\boldsymbol{\mathsf A}}} 
\def \CC {{\boldsymbol{\mathsf C}}} 
\def \UU {{\boldsymbol{\mathsf U}}} 
\def \Rset {{\mathbb R}}         %mathsets

\def \Zset {{\mathbb Z}}
\def \Nset {{\mathbb N}}
\def \Qset {{\mathbb Q}}

\def\BA{\mathbb A}
\def \FF {{\mathcal{F}}}

\def \TT {{\mathcal{T}}}

\def \ZZ {{\mathcal{Z}}}

\def \al {\alpha}
\def \be {\beta}

\def \La {\Lambda}

\def \ga {\gamma}

\def \ep {\varepsilon}
\def \lra {\longrightarrow}

\def \hra {\hookrightarrow}

\def \rcor {\rangle}
\def \lcor {\langle}
\def \ol {\overline}
\def \wt {\widetilde}

\def \Id { {\mathrm{Id}} }

\def \mm  {\mathfrak{m}}

\DeclareMathOperator \kk {\Bbbk}
\DeclareMathOperator \tr { {\mathrm{tr}} }
\DeclareMathOperator \Tr { {\mathrm{Tr}} }
\DeclareMathOperator \sta { {\mathrm{st}} }
\DeclareMathOperator \reg { {\mathrm{reg}} }
\DeclareMathOperator \red{ {\mathrm{red}} }
\DeclareMathOperator \cc { {\mathrm{c}}}
\DeclareMathOperator \nc { {\mathrm{nc}}}
\DeclareMathOperator \Span { {\mathrm{Span}} }

\DeclareMathOperator \charr { {\mathrm{char}} }

\DeclareMathOperator \ord { {\mathrm{ord}} }

\DeclareMathOperator \diag { {\mathrm{diag}} }
\DeclareMathOperator \Ker { {\mathrm{Ker}} }

\DeclareMathOperator \End { {\mathrm{End}} }

\DeclareMathOperator \MaxSpec { {\mathrm{MaxSpec}}}
\DeclareMathOperator \Ray{\mathrm{Ray}}
\renewcommand \Im { {\mathrm{Im}} }

\newcommand \supp { {\mathrm{supp}} }
\newcommand \Fract { {\mathrm{Fract}} }

\def \B  {{\widetilde{B}}}
\newcommand \ex {{\bf{ex}}}
\newcommand \inv {{\bf{inv}}}

\def\uL{\underline{\Phi}}
\def\CL{\mathcal L}
\def\Fr{\mathrm{Fr}}
\def\Newt{\mathrm{Newt}}
\begin{document}

\title[Cluster Algebras and Cayley--Hamilton Algebras]{Root of Unity Quantum Cluster Algebras and Cayley--Hamilton Algebras}
\author[S. Huang]{Shengnan Huang}
\thanks{The research of S.H. and M.Y. has been supported by NSF grants DMS-1901830 and DMS-2131243. The research of T.L. has been supported by NSF grant DMS-1811114.}
\address{
Department of Mathematics, 360 Huntington Avenue, Northeastern University, Boston, MA 02115, USA}
\email{huang.shengn@northeastern.edu}
\author[T. T. Q. L\^{e}]{Thang T. Q. L\^{e}}
\address{
School of Mathematics, 686 Cherry Street, Georgia Tech, Atlanta, GA 30332, USA}
\email{letu@math.gatech.edu}
\author[M. Yakimov]{Milen Yakimov}
\address{
Department of Mathematics, 360 Huntington Avenue, Northeastern University, Boston, MA 02115, USA}
\email{m.yakimov@northeastern.edu}
\keywords{Quantum cluster algebras at roots of unity, algebras with trace, Cayley--Hamilton algebras, maximal orders}
\subjclass[2010]{Primary: 13F60, Secondary: 16G30, 17B37, 14A22}
\begin{abstract} 
We prove that large classes of algebras in the framework of root of unity quantum cluster algebras
have the structures of maximal orders in central simple algebras and Cayley--Hamilton algebras in the sense of Procesi.
We show that every root of unity upper quantum cluster algebra is a maximal order and obtain an explicit 
formula for its reduced trace. Under mild assumptions, inside each such algebra we construct a 
canonical central subalgebra isomorphic to the underlying upper cluster algebra, 
such that the pair is a Cayley--Hamilton algebra; its fully Azumaya locus is shown 
to contain a copy of the underlying cluster $\mathcal{A}$-variety. Both results are proved in the wider generality of 
intersections of mixed quantum tori over subcollections of seeds. Furthermore, we prove that 
all monomial subalgebras of root of unity quantum tori are Cayley--Hamilton algebras and 
classify those ones that are maximal orders. Arbitrary intersections of those over subsets 
of seeds are also proved to be Cayley--Hamilton algebras. Previous approaches to 
constructing maximal orders relied on filtration and homological methods. We use new 
methods based on cluster algebras.
\end{abstract}
\maketitle
%%%%%%%%%%%%%%%%%%%%%%%%%%%%%%%%%%%%%%%%%%%%%%%%%%%%%%%%%%%%%%%%%%%%%%%%%%%
\sectionnew{Introduction}
\lb{Intro}
\subsection{Cayley--Hamilton algebras} Cayley--Hamilton algebras, defined by Procesi \cite{P} in the late 80s, provide a vast generalization of maximal orders in central simple algebras.
Such an algebra of degree $d$ (a positive integer) is a triple $(R,C, \tr)$ consisting of a $k$-algebra $R$ of a commutative ring $k$, a central subalgebra $C$ and a trace function $\tr$ satisfying the 
$d$-th Cayley--Hamilton identity, which means that each element of $R$ satisfies its $d$-th characteristic polynomial. Every maximal order in a central simple algebra of PI degree $d$ whose 
center has characteristic $p \notin [1,d]$ is a Cayley--Hamilton algebra of degree $d$ with respect to its full center and reduced trace, cf. Sect. \ref{traces}.

There are a number of results on Cayley--Hamilton algebras that provide a general approach to their representation theory. Using invariant theory, for each character 
$\chi$ of $C$, Procesi \cite{P, DP} constructed a semisimple representation of $R$ whose direct summands exhaust all irreducible representations of $R$ with central character $\chi$. 
In \cite{BY} it was proved that there is a close relationship between the discriminant ideals of a Cayley--Hamilton algebras and its irreducible representations,
namely the zero loci of the discriminant ideals record the sum of the squares of the irreducible representations with a given central character $\chi$.

There has been a great interest in proving that quantum algebras at roots of unity that appear in Lie theory and topology posses structures of maximal orders or more generally 
Cayley--Hamilton algebras with the aim of classifying their irreducible representations. 

In Lie theory, De Concini, Kac and Procesi \cite{DK,DP} proved that all big quantum groups at roots of unity and all 
quantum Schubert cells are maximal orders in central simple algebras. Their centers are singular which presents difficulties for  
applications to representation theory. 
To overcome those, it was proved in \cite{DK,DP} that each algebra in the two classes 
possesses a non-singular subalgebra and an appropriate trace function that makes the pair a Cayley--Hamilton algebra.
Moreover, De Concini and Lyubashenko \cite{DL} proved that the quantum function algebras at roots of unity of all complex simple Lie groups 
are maximal orders. Stafford proved that the 3 and 4-dimensional Sklyanin algebras corresponding to elliptic curves with 
finite order automorphisms are maximal orders \cite{S}. 
 
In topology, it was proved that the stated skein algebras at roots of unity of surfaces with nontrivial boundary \cite{LY} 
are maximal orders. In \cite{Paprocki} it was proved that in those cases the Muller skein algebras \cite{Muller} at roots of unity 
are also maximal orders. 

All of the above results on maximal orders are proved by filtration arguments, homological methods or by methods relying on 
unique factorization domains. The results on Cayley--Hamilton structures on quantum groups at roots of unity
construct Cayley--Hamilton algebras $(R,C, \tr)$ of a very special nature having the property that $R$ is a free $C$-module. 

\subsection{From root of unity upper quantum cluster algebras to Cayley--Hamilton algebras and maximal orders} In this paper we prove that 
all root of unity upper quantum cluster algebras are maximal orders in central simple algebras (maximal orders for short).
This is a very general setting that unifies quantum algebras from Lie theory, topology and mathematical physics.

Furthermore, under mild assumptions on the order of the root of unity, we prove that each 
root of unity upper quantum cluster algebra possesses a canonical central subalgebra isomorphic to 
the corresponding upper cluster algebra and a trace function making the triple into a Cayley--Hamilton algebra. Unlike the situation 
in the work of De Concini, Kac and Procesi, for our triples $(R,C, \tr)$, $R$ is rarely a free module over $C$. 

Cluster algebras have been an area of intense research in the last 20 years due to their deep relations to many areas of mathematics and mathematical physics, 
see \cite{FWZ,GSV,M}. In this paper we use the methods of cluster algebras to construct broad classes of maximal orders and Cayley--Hamilton algebras.
Previously, structures of maximal orders were constructed using filtration arguments \cite{MR} and homological methods \cite{S}. However, those do not apply to the classes that we construct.
In the broad generality that is considered in the paper, root of unity upper quantum cluster algebras do not posses any filtrations 
that can be used to prove that they are maximal orders. They do not satisfy the homological assumptions in \cite{S}; for instance, they are 
Auslander regular only in special cases. 

\subsection{Statements of main results} Fix an integral domain $\kk$ of characteristic 0. 
Let $\ell$ be a positive integer, $\ep^{1/2}$ be a primitive $\ell$-th root of unity 
in the algebraic closure of the fraction field of $\kk$, and
\[
\Abbe := \kk[\ep^{1/2}].
\]
To a compatible pair, consisting of an exchange matrix $\B$ and a root of unity toric frame $M_\ep$,
one associates \cite{NTY} the root of unity upper quantum cluster algebra $\UU_\ep(M_\ep,  \B, \inv)$, where $\inv$ is a subset of the set of 
frozen variables denoting those that are inverted, see Sect. \ref{2.3} for details. It is the $\Abbe$-algebra given by the intersection
\[
\UU_\ep(M_\ep, \B, \inv) := \hspace{-15pt} \bigcap_{\mbox{all seeds} \; (M'_\ep, \B')} \hspace{-15pt} \TT_\ep(M'_\ep)_{\geq}
\]
where $\TT_\ep(M'_\ep)_{\geq}$ is the mixed quantum torus associated to the toric frame $M'_\ep$ which is the $\Abbe$-subalgebra of the root of unity quantum 
torus $\TT_\ep(M'_\ep)$ corresponding to $M'_\ep$ obtained by not inverting the frozen variables that do not lie in $\inv$, cf. \eqref{TM}--\eqref{TMgeq}.
Here the term {\em{mixed}} refers to the fact that $\TT_\ep(M'_\ep)_{\geq}$ is a mixture of a quantum torus and a quantum plane.

For a subset of seeds $\Theta$, consider the algebra 
\[
\UU_\ep(M_\ep, \B, \inv, \Theta) := \hspace{-15pt} \bigcap_{(M'_\ep, \B') \in \Theta} \hspace{-15pt} \TT_\ep(M'_\ep)_{\geq}.
\] 
Our results concern these more general algebras, rather than just root of unity upper quantum cluster algebras $\UU_\ep(M_\ep, \B, \inv)$.
This adds an extra degree of flexibility for applications, since potentially there could be interesting algebras that are isomorphic to algebras of the 
form $\UU_\ep(M_\ep, \B, \inv, \Theta)$ but are not isomorphic to a root of unity upper quantum cluster algebra. 
\medskip

\noindent
{\bf{Theorem A.}}
{\em{Let $\Theta$ be an arbitrary subset of seeds of $\UU_\ep(M_\ep, \B, \inv)$. The following hold:
\begin{enumerate}
\item The algebra $\UU_\ep(M_\ep, \B, \inv, \Theta)$ is a maximal order in a central simple algebra. Its reduced trace is the 
restriction of any of the reduced traces of the quantum tori $\TT_\ep(M'_\ep)$, 
explicitly given by \coref{red-tr-quant-tor}, for the seeds $(M'_\ep, \B')$ in $\Theta$.
\item If the base ring $\kk$ is an algebraically closed field containing $\Qset(\ep^{1/2})$, then the union
\[
\bigcup_{(M'_\ep, \B') \in \Theta} \MaxSpec \big( \ZZ(\UU_\ep(M_\ep, \B, \inv, \Theta))[ M'_\ep(e_i)^{-\ell}, 1 \leq i \leq N] \big)
\]
inside $\MaxSpec \ZZ(\UU_\ep(M_\ep, \B, \inv, \Theta))$
is in the Azumaya locus of $\UU_\ep(M_\ep, \B, \inv, \Theta)$.
\end{enumerate}
}}
\medskip
We follow the conventions of \cite{MR} for orders in central simple algebras 
with equivalent definitions given in \cite[Ch. 3, \S 1.2]{MR} and \cite[Ch. 5 \S 3.6]{MR}, and equivalence proved in 
\cite[Proposition 5.3.10]{MR}. We refer the reader to Sect. \ref{sec.maxorder} for all necessary background on maximal orders.

For a PI algebra, its Azumya locus is an important representations theoretic object, classifying its 
irreducible representations of maximal dimension, cf. \deref{fullyAzumaya}. 

For a mixed quantum torus $\TT_\ep(M_\ep)_\geq$, let $\TT_\ep^\ell(M_\ep)_\geq$ be its central $\Abbe$-subalgebra generated 
by the $\ell$-th powers of the standard generators of $\TT_\ep(M_\ep)_\geq$. It is a mixture of a Laurent polynomial ring and polynomial ring. 
$\TT_\ep(M'_\ep)_{\geq}$ is free over $\TT_\ep^\ell(M'_\ep)_{\geq}$ which gives rise to the regular trace function 
\[
\tr_{\reg}^{\TT_\ep(M'_\ep)_{\geq}} : \TT_\ep(M'_\ep)_{\geq} \to \TT_\ep^\ell(M'_\ep)_{\geq},
\]
explicitly described in \leref{T-Tl}. Define 
\[
\CC\UU_\ep(M_\ep, \B, \inv, \Theta) := \hspace{-15pt} \bigcap_{ \ \ (M'_\ep, \B') \in \Theta} \hspace{-15pt} \TT_\ep(M'_\ep)^\ell_{\geq}.
\]

If $\ell$ is odd and coprime to the skew-symmetrizing integers for the principal part of $\B$
and $\Theta$ is the set of all roots, then $\CC\UU_\ep(M_\ep, \B, \inv, \Theta)$ is isomorphic to 
the corresponding upper cluster algebra defined over $\Abbe$:
\[
\CC\UU_\ep(M_\ep, \B, \inv) \cong \UU(\B, \inv)_{\Abe},
\]
see \prref{cent-U} below. In \prref{cent-U-Theta} an extension of this fact to the algebras $\CC\UU_\ep(M_\ep, \B, \inv, \Theta)$ 
and their classical counterparts in the setting of upper cluster algebras is proved. 
\medskip

\noindent
{\bf{Theorem B.}}
{\em{Assume that the order $\ell$ of the root of unity $\ep^{1/2}$ is odd and coprime to the skew-symmetrizing integers for the principal part of $\B$. 
Let $\Theta$ be a connected set of vertices of the exchange graph of $\UU_\ep(M_\ep, \B, \inv)$. The following hold:
\begin{enumerate}
\item
For every pair of seeds $(M'_\ep, \B'), (M''_\ep, \B'') \in \Theta$, 
\[
\tr_{\reg}^{\TT_\ep(M'_\ep)_{\geq}}\big{|}_{\UU_\ep(M_\ep, \B, \inv, \Theta)} = \tr_{\reg}^{\TT_\ep(M''_\ep)_{\geq}}\big{|}_{\UU_\ep(M_\ep, \B, \inv, \Theta)}.
\]
Denote by $\tr_{\reg}$ this restriction map coming from an arbitrary seed in $\Theta$. 
\item $\tr_{\reg} ( \UU_\ep(M_\ep, \B, \inv, \Theta)) \subseteq \CC\UU_\ep(M_\ep, \B, \inv, \Theta)$.
\item The triple 
\[
(\UU_\ep(M_\ep, \B, \inv, \Theta), \CC\UU_\ep(M_\ep, \B, \inv, \Theta), \tr_{\reg}) 
\]
is a Cayley--Hamilton algebra of degree equal to $\ell^N$, where $N$ is the number of cluster variables.
\item Assume that the base ring $\kk$ is a field extension of the cyclotomic field $\Qset(\ep^{1/2})$.
Then $\UU_\ep(M_\ep, \B, \inv, \Theta)$ is a finitely generated $\kk$-algebra if and only 
$\CC \UU_\ep(M_\ep, \B, \inv, \Theta)$ is a finitely generated $\kk$-algebra and 
$\UU_\ep(M_\ep, \B, \inv, \Theta)$ is a finitely generated module over $\CC \UU_\ep(M_\ep, \B, \inv, \Theta)$.
\item If the base ring $\kk$ is an algebraically closed field containing $\Qset(\ep^{1/2})$, then the union
\begin{equation}
\label{full-Azum}
\bigcup_{(M'_\ep, \B') \in \Theta} \MaxSpec \big( \CC\UU_\ep(M_\ep, \B, \inv, \Theta)[ M'_\ep(e_i)^{-\ell}, 1 \leq i \leq N] \big)
\end{equation}
inside $\MaxSpec \CC \UU_\ep(M_\ep, \B, \inv, \Theta)$
is in the fully Azumaya locus of the algebra $\UU_\ep(M_\ep, \B, \inv, \Theta)$
with respect to the central subalgebra $\CC \UU_\ep(M_\ep, \B, \inv, \Theta)$,
see \deref{fullyAzumaya}.
\end{enumerate}
}}
\medskip

Part (5) of the Theorem B provides yet another link between noncommutative ring theory and cluster algebras:
when $\Theta$ is the set of all seeds, \eqref{full-Azum} is precisely the cluster $\mathcal{A}$-variety 
of the associated cluster algebra defined over $\kk$, see \reref{cluster-A-variety}.
In other words, the fully Azumaya locus of the root of unity quantum upper cluster algebra
$\UU_\ep(M_\ep, \B, \inv)$ contains a copy of the associated cluster $\mathcal{A}$-variety over $\kk$.
Undoubtably, cluster algebra theory will play a role in classifying all irreducible representations of the algebras that are isomorphic to algebras 
of the form $\UU_\ep(M_\ep, \B, \inv, \Theta)$.

Up to date there are practically no general results that transfer properties between upper cluster algebras and their quantum and root of unity 
counterpart. Part (4) of Theorem B proves such a result. It implies the following fact:
\medskip

\noindent
{\bf{Corollary.}} {\em{If the root of unity quantum cluster algebra $\UU_\ep(M_\ep, \B, \inv)$ is a finitely generated algebra over a field extension of $\Qset(\ep^{1/2})$, 
then the upper cluster algebra $\UU_\ep(M_\ep, \B)$ has the same property.}}
\medskip

In Theorems A and B we intersect mixed quantum tori because those are the natural objects that appear in cluster algebra theory. 
But our methods allow us to deal with quite more general classes of algebras.
For a seed $\Sigma = (M_\ep, \B)$ and a submonoid $\Phi$ of $\Zset^N$, denote by $\BA_\Sigma(\Phi)$ the $\Abbe$-subalgebra 
of $\TT_\ep(M_\ep)$ generated by the monomials with exponents in $\Phi$. We call those monomial subalgebras.
For example, the mixed quantum torus $\TT_\ep(M_\ep)_\geq$ is a special case of a subalgebra of $\TT_\ep(M_\ep)$ of this type.
\medskip

\noindent
{\bf{Theorem C.}}
{\em{
\begin{enumerate}
\item For every monomial subalgebra $\BA_\Sigma(\Phi)$ of $\TT_{\ep}(M_\ep)$, the triple
\[
\big( \BA_\Sigma(\Phi), \BA_\Sigma(\Phi) \cap \ZZ(\TT_{\ep}(M_\ep)), \tr_{\red}^{\TT_{\ep}(M_\ep)} \big)
\]
is a Cayley--Hamilton algebra of degree given in \prref{monoid-CH-2}. 
\item A monomial subalgebra $\BA_\Sigma(\Phi)$ is a maximal order if and only if the submonoid $\Phi$ of $\Zset^N$ is integrally convex and integrally closed,
see Sect. \ref{monoid-max-ord} for definitions. 
\end{enumerate}

For the remaining part of the theorem, let $\Theta$ be a subset of seeds of $\UU_\ep(M_\ep, \B, \inv)$ and 
$\BA_\Sigma(\Phi_\Sigma)$ be a monomial subalgebra of $\TT_\ep(M'_\ep)$ for each seed $\Sigma = (M'_\ep, \B') \in \Theta$.
Consider the $\Abbe$-algebra
\[
\BA := \bigcap_{\Sigma \in \Theta} \BA_\Sigma(\Phi_\Sigma). 
\]
The following hold:
\begin{enumerate}
\item[(3)] For every pair of seeds $(M'_\ep, \B'), (M''_\ep, \B'') \in \Theta$, 
\[
\tr_{\red}^{\TT_\ep(M'_\ep)}\big{|}_{\BA} = \tr_{\red}^{\TT_\ep(M''_\ep)}\big{|}_{\BA}.
\]
Denote by $\tr_{\red}$ this restriction map coming from an arbitrary seed in $\Theta$. 
\item[(4)] $\tr_{\red} ( \BA ) \subseteq \cap_{\Sigma \in \Theta} \ZZ(\BA_\Sigma(\Phi_\Sigma))$.
\item[(5)] The triple 
\[
(\BA, \cap_{\Sigma \in \Theta} \ZZ(\BA_\Sigma(\Phi_\Sigma)), \tr_{\red}) 
\]
is a Cayley--Hamilton algebra of degree given in \thref{CH-U2}.  
\end{enumerate}
}}
\medskip

Finally, we note that when algebras are proved to be maximal orders using filtration arguments \cite{MR} and homological methods \cite{S}, this 
does not give information on their reduced traces. Unlike those methods, in Theorems A-C we obtain explicit formulas for the 
corresponding trace maps. All results in the paper are proved integrally over $\kk[\ep^{1/2}]$
for an arbitrary integral domain $\kk$ of characteristic 0,
except part (2) of Theorem A and parts (4)-(5) of Theorem B.

The paper is organized as follows. Sect. 2 contains background material on Cayley--Hamilton algebras and cluster algebras.
Sect. 3 sets up the framework for root of unity upper quantum cluster algebras, the algebras of the form $\UU_\ep(M_\ep, \B, \inv, \Theta)$, 
and their central subalgebras.
The proofs of parts (1)-(4) of Theorem B is given in Sect. 4. Sect. 5 contains background material on maximal orders and reduced traces and
proves a general theorem on intersections of maximal orders and Caylay--Hamilton algebras. In Sect. 6 we prove Theorem A and 
part (5) of Theorem B.
Sect. 7 illustrates Theorems A and B with a root of unity upper quantum cluster algebra
without frozen variables; examples that were treated before in Lie theory and topology always have 
sufficient number of frozen variables which makes them more tractable. Sect. 8 contains our results on 
monomial algebras and their cluster theoretic intersections; it proves Theorem C.

For additional background material on maximal orders, representation theory of PI algebras, and Cayley--Hamilton algebras we refer the reader to 
\cite{MR,Re}, \cite{BG,MR} and \cite{DP}, respectively.
\medskip

\noindent
{\bf{Notation:}} The center of a ring $R$ will be denoted by $\ZZ(R)$. For a commutative ring $A$, we will denote by $M_n(A)$ the ring of $n \times n$ matrices with entries in $R$ and by
\[
\Tr : M_n(A) \to A
\]
the standard matrix trace. We will denote by $\Nset$ the set of non-negative integers and by $\Zset_+$ the set of positive integers. 
\medskip

\noindent
{\bf{Acknowledgements.}} We are grateful to Harm Derksen, Visu Makam and Alex Martsinkovsky for valuable discussions.  
%%%%%%%%%%%%%%%%%%%%%%%%%%%%%%%%%%%%%%%%%%%%%%%%%%%%%%%%%%%%%%%%%%%%%%%%%%
\sectionnew{Preliminaries on Cayley--Hamilton algebras and cluster algebras}
\lb{background}
This section contains preliminaries on Cayley--Hamilton algebras, cluster algebras of geometric type,
root of unity quantum cluster algebras and exchange graphs that will be used later in the paper.
\subsection{Cayley--Hamilton algebras}
For $1 \leq i \leq d$, denote by $\sigma_i$ the $i$-th elementary symmetric function in the indeterminates 
$\lambda_1, \lambda_2, \dots, \lambda_d$ and by $\psi_i:=\lambda_1^i+\lambda_2^i+\cdots +\lambda_d^i$ the Newton power sum function.
It is well known that there exists a unique set of polynomials 
\[
p_i(x_1, x_2, \dots, x_i)\in \mathbb{Z}[(i!)^{-1}][x_1, x_2, \dots, x_i]
\] 
such that
\[
\sigma_i=p_i(\psi_1, \psi_2, \dots, \psi_i), \quad \forall \; 1\leq i \leq d.
\]

\bde{alg-tr}
An algebra with trace is a triple $(R, C, \mathrm{tr})$, where $R$ is a $\kk$-algebra with $\kk$ being a commutative ring, 
$C$ is a central subalgebra of $R$ and $\tr: R \to C$ is a $C$-linear map such that
\[
\tr(ab) = tr(b a) \quad \mbox{for all} \quad a, b \in R.
\]
\ede
Fix a positive integer $d$ and assume that $i$ is not a zero divisor of $R$ for $1 \leq i \leq d$. The $d$-th characteristic polynomial $\chi_{d,a}(t) \in C[(d!)^{-1}][t]$ 
of an element $a \in R$ is defined to be
\[
\chi_{d, a}(t):= t^d-c_1(a)t^{d-1}+\cdots +(-1)^d c_d(a),
\]
where $c_i(a):=p_i\big(\tr(a), \tr(a^2), \dots, \tr(a^i)\big)$.

\bde{CHring} A Cayley--Hamilton algebra of degree $d \in \Zset_+$ is a $\kk$-algebra with trace $(R, C, \tr)$
over a commutative ring $\kk$ 
such that $i$ is a not a zero divisor of $R$ for $1 \leq i \leq d$ and 
\begin{enumerate}
\item for all $a\in R$, $\chi_{d,a}(a)=0$,
\item $\tr(1)=d$.
\end{enumerate}
\ede

\bex{max-ord}
Every maximal order, see Sect. \ref{sec.maxorder} below, of PI degree $d$ whose center has characteristic $p \notin [1,d]$ 
is a Cayley--Hamilton algebra of degree $d$
with respect to the reduced trace, see Sect. \ref{traces} for details.
\eex

\ble{restrict-CH}
If $(R, C, \tr)$ is a Cayley--Hamilton algebra of degree $d$ and $R'$, $C'$ are subalgebras of $R$, $R' \cap C$, respectively such that $\tr(R') \subseteq C'$, 
then $(R', C', \tr|_{R'})$ is also a Cayley--Hamilton algebra of degree $d$. 
\ele
\begin{proof} Since $\tr : R \to C$ is $C$-linear and $C' \subseteq R' \cap C$, $\tr|_{R'} : R' \to C'$ is $C'$-linear. The conditions 
(1)-(2) in \deref{CHring} for $\tr$ imply their validity for $\tr|_{R'}$.
\end{proof}
\bde{fullyAzumaya} Let $R$ be a prime affine algebra over an algebraically closed field $\kk$, 
which is a finitely generated module over a central $\kk$-subalgebra $C$. 
\begin{enumerate}
\item A point $\mm \in \MaxSpec Z(R)$ is in the {\em{Azumaya locus}} of $R$ if $R_\mm$ is an Azumaya algebra over
$Z_\mm$. This is equivalent to saying that $R$ has an irreducible module, annihilated by $\mm$, of maximal dimension among the irreducible $R$-modules
(which equals the PI degree of $R$); such a representation is automatically unique (see \cite[Theorem III.1.6]{BG}).

\item \cite{BrGo} A maximal ideal $\mm$ of $C$ is in the fully Azumaya locus of $R$ with respect to $C$ if all irreducible representations of $R$ that are annihilated by $\mm R$ have maximal 
dimension  

This can be formulated in the following equivalent way. Denote by 
\[
\psi: \MaxSpec \ZZ(R) \to \MaxSpec C 
\]
the map induced by the inclusion $C \hra Z(R)$. A maximal 
ideal $\mm \in \MaxSpec C$ is in the fully Azumaya locus of $R$ with respect to $C$ if all preimages in $\psi(\mm)$ are in 
the Azumaya locus of $R$,  cf. \cite[Sect. III.1.7]{BG}.  
\end{enumerate}
\ede

\subsection{Maximal orders} \label{sec.maxorder}
We recall the notion of a maximal order, following \cite{MR}.
A {\em central simple algebra } $S$ over a commutative field $Q$ is a  finite dimensional simple $Q$-algebra whose center is $Q$. 
An {\em order} in a central simple algebra $S$ is a subring $R\subseteq S$  such that every $x\in S$ can be presented as $x = a b^{-1} = c^{-1} d$ with $a,b,c,d\in R$,
see \cite[Ch. 3, \S 1.2]{MR}. Two orders $R, R'$ of $S$ are {\em equivalent} if there are units $a,b,a',b'\in S$ such that $aRb \subseteq R'$ and $a'R'b' \subseteq R$. 
A {\em maximal order} is an order which is maximal (with respect to inclusions) in its equivalence class. 
Actually the notion of an order is more generally defined for quotient rings, see \cite[Ch. 3, \S 1.2]{MR}, but 
in the case of central simple algebras, there is a simpler equivalent definition. The following is essentially given in \cite[Sect. 6]{DP}.

\bpr{r.maxorder} 
Assume $S$ is a central simple algebra over a commutative field $Q$ and $R \subseteq S$ is a subring with center $\ZZ(R)$.
\begin{enumerate}
\item $R$ is an order of $S$ if and only if 
\begin{enumerate}
\item[(i)] $\ZZ(R)$ is an integral domain whose field of fractions is $Q$ and 
\item[(ii)] $RQ=S$.
\end{enumerate}
\item If $R$ is an order, then it is a maximal order if and only if 
\begin{itemize}
\item[(iii)] for any ring $R'$ with 
\begin{equation}
  R \subseteq R' \subseteq \frac{1}{z} R := \{ x z^{-1} \mid x\in R\}, \label{r.maxord2}
\end{equation}
where $z\in \ZZ(R)$ is a non-zero element, we have $R'=R$.
\end{itemize}
\end{enumerate} 
\epr

\begin{proof} (1) Suppose $R$ is an order in $S$. By \cite[Proposition 5.3.10]{MR} the ring $R$ is another type or order, called $\ZZ(R)$-order there,
see \cite[Ch. 5 \S 3.6]{MR}, and (i) and (ii) are part of the definition of a $\ZZ(R)$-order.
Thus we have (i) and (ii).
 
 The converse follows obviously from the definition of order. 
 
 (2) Assume that $R$ is a maximal order and $R'$ satisfies \eqref{r.maxord2}. 
 As $R \subseteq R'$, by \cite[Corollary 3.1.6]{MR}, the ring $R'$ is also an order. 
 Eq. \eqref{r.maxord2} implies $R'$ is equivalent to $R$, and maximality means $R'=R$.
 
 Let us prove the converse. Assume (iii). Let $R'$ be an order equivalent to $R$. By Proposition \cite[Proposition 5.3.8(iii)]{MR}, there is a non-zero $z\in \ZZ(R)$ such that $R' \subseteq \frac 1z  R$. If in addition $R \subseteq R'$, then by (iii) we have $R'=R$. This shows that $R$ is a maximal order.
\end{proof}
\bre{second-def-order}
In some texts like Reiner's book \cite{Re}, 
the notion of a maximal order $R$ in a central simple algebra $Q$ is stronger and requires the center $\ZZ(R)$ to be Noetherian.
In the definition of $R$ being a $\ZZ(R)$-order in $Q$ one adds the condition that $R$ be a finitely generated $\ZZ(R)$-module, 
\cite[Sect. 8]{Re}; this notion is called {\em{classical order}} in \cite[Ch. 5 \S 3.5]{MR}.
However, under the assumption that $\ZZ(R)$ is Noetherian, this condition follows from 
condition (ii) in \prref{r.maxorder}(1), see \cite[Proposition 5.3.14]{MR}. 
\ere
\subsection{Cluster algebras of geometric type}
Let $\kk$ be an integral domain of characteristic 0 and $F$ be its fraction field.

Let $N$ be a positive integer, $\ex\subseteq [1,N]$ (set of {\em{exchangeable indices}}), $\inv \subseteq [1,N] \backslash \ex$ (set of {\em{frozen variables}} that will be {\em{inverted}}) and 
$\FF$ be a purely transcendental extension of $F$ of degree $N$. A {\em{seed}} is a pair $(\wt{\mathbf{x}}, \B)$ such that 
\begin{enumerate}
\item $\wt{\mathbf{x}}=\{x_1, \dots, x_N \}$ is a transcendence basis of $\FF$ over $F$;
\item $\B = (b_{ij}) \in M_{N\times \ex}(\Zset)$ ({\em{exchange matrix}}) and its $\ex \times \ex$ submatrix $B$ (the principal part of $\B$)
is skew-symmetrizable. More precisely, $D B$ is skew-symmetric for a matrix $D = \diag(d_j, j \in \ex)$ with $d_j \in \Zset_+$. 
\end{enumerate}
The mutation of the seed $(\wt{\mathbf{x}}, \B)$ in the direction of $k \in \ex$ is defined to be the seed $\mu_k(\wt{\mathbf{x}}, \B):= (\wt{\mathbf{x}}', \B')$, where  
\begin{equation}
\label{classical-mutation}
\wt{\mathbf{x}}' = \{x_k'\} \cup \wt{\mathbf{x}} \backslash \{x_k\} 
\quad
\mbox{and}
\quad
x_k x_k' := \prod_{b_{ik}>0}x_i^{b_{ik}} + \prod_{b_{ik}<0}x_i^{-b_{ik}}
\end{equation}
and 
\begin{equation}
\label{mukB}
\mu_k(\B)= (b_{ij}'):=
\begin{cases}
-b_{ij}, & \text{if } i = k \text{ or } j=k \\
b_{ij} + \frac{|b_{ik}|b_{kj} + b_{ik}|b_{kj}|}{2}, & \text{otherwise.}
\end{cases}
\end{equation}
The principal part of $\B'$ is skew-symmetrized by the same matrix $D$. Mutation is involutive: $\mu^2_k(\wt{\mathbf{x}}, \B) = (\wt{\mathbf{x}}, \B)$. 
Two seeds $(\wt{\mathbf{x}}, \B)$ and $(\wt{\mathbf{x}}'', \B'')$ are {\em{mutation-equivalent}}, 
denoted $(\wt{\mathbf{x}}, \B) \sim (\wt{\mathbf{x}}'', \B'')$, if one is obtained from the other by a sequence of mutations.  

\bde{cluster algebra} \cite{FZ1,BFZ}
\hfill 
\begin{enumerate}
\item The {\em{cluster algebra}} $\AA (\B, \inv)_{\kk}$ is the unital $\kk$-subalgebra of $\FF$ generated by the cluster variables in the seeds 
$(\wt{\mathbf{x}}'', \B'') \sim (\wt{\mathbf{x}}, \B)$ and by $x_i^{-1}$ for $i \in \inv$. 
\item The upper cluster algebra $\UU(\B, \inv)_{\kk}$ is the intersection of all mixed polynomial/Laurent polynomial rings
\[
\UU(\B, \inv)_{\kk} := \bigcap_{(\wt{\mathbf{x}}'', \B'') \sim (\wt{\mathbf{x}}, \B)} \kk[(x''_i)^{\pm 1}, x''_j; i \in \ex \sqcup \inv, j \notin \ex \sqcup \inv].
\]
\end{enumerate}
\ede
It is clear that
\begin{equation}
\label{base-change}
\AA(\B, \inv)_{\kk} \cong \kk \otimes_\Zset \AA(\B, \inv)_\Zset
\end{equation}
but we are not aware of a similar fact for $\UU(\B, \inv)_{\kk}$ unless we are in the case when $\kk$ is a finite extension of $\Zset$. 

{\em{When the base ring $\kk$ is clear from the discussion, we will use the notations $\AA (\B, \inv)$ and $\UU(\B, \inv)$ 
for brevity, but $\kk$ will be an arbitrary integral domain of characteristic 0 and not just $\Zset$.}}

By the Laurent Phenomenon Theorem of Fomin and Zelevinsky \cite{FZlp}, 
\[
\AA (\B, \inv) \subseteq \UU(\B, \inv). 
\] 

\subsection{Root of unity quantum cluster algebras}
\label{2.3}
We follow the framework of Berenstein--Zelevinsky for quantum cluster algebras \cite{BZ}, adapted to the root of unity case in \cite{NTY}. The 
algebras that are considered are in general not specializations of quantum cluster algebras. This treatment can be viewed as defining quantum cluster 
$\mathcal{A}$-varieties at roots of unity. Quantum cluster $\mathcal{X}$-varieties at roots of unity were defined and studied by 
Fock and Goncharov in \cite{FG}.

As in the previous subsection, $\kk$ will denote an integral domain of characteristic 0 and $F$ its fraction field.
Fix a positive integer $\ell$ and a primitive $\ell$-th root of unity, $\ep^{1/2}$ in the algebraic closure of $F$. Recall from the introduction that
\[
\Abbe := \kk[\ep^{1/2}].
\]
Let 
\[
\La: \Zset^N \times \Zset^N \to \Zset/\ell := \Zset/ (\ell \Zset)
\]
be a skew-symmetric bicharacter. The root of unity (based) quantum torus $\TT_{\ep}(\La)$ 
is the $\Abbe$-algebra with basis $\{ \hspace{1pt} X^f \hspace{1pt} | \hspace{1pt} f \in \Zset^N \hspace{1pt} \}$ and product
\begin{equation}
\label{prod-quant-tor}
X^f X^g = \ep^{\La(f,g)/2} X^{f+g} \quad \mbox{for all} \quad f,g \in \Zset^N.
\end{equation}
Obviously $\TT_\ep(\La)$ is a domain. A {\em root of unity toric frame} $M_\ep$ of a division algebra $\FF_\ep$ over $\Qset(\ep^{1/2})$ is a map $M_\ep : \Zset^N \to \FF_\ep$ 
for which there exists a bicharacter $\La$ as above with the properties:
\begin{enumerate}
\item There is an $\Abbe$-algebra embedding $\phi: \TT_\ep(\La) \hra \FF_\ep$ given by $\phi(X^f)=M_\ep(f)$, $\forall f\in \Zset^N$.
\item $\FF_\ep \simeq \Fract \left( \phi( \TT_\ep(\La) \right) )$.
\end{enumerate}

The matrix $\La \in M_N(\Zset/\ell)$ (called {\em{matrix of the frame}} $M_\ep$) is uniquely reconstructed from the root of unity toric frame $M_\ep$. Denote
quantum torus 
\begin{equation}
\label{TM}
\TT_\ep(M_\ep) = \phi(\TT_\ep(\La)) \subseteq \FF_\ep.
\end{equation}
Define the {\em{mixed quantum tori}}
\begin{equation}
\label{TMgeq}
\TT_\ep(\La)_{\geq} \quad \mbox{and} \quad \TT_\ep(M_\ep)_{\geq}
\end{equation}
inside $\TT_\ep(\La)$ and  $\TT_\ep(M_\ep)$ with 
$\Abbe$-bases consisting of the elements $X^f$ and $M_\ep(f)$ for those $f = (f_1, \ldots f_N) \in \Zset^N$ such that $f_j \geq 0$ for $j \notin \ex \sqcup \inv$.
The restriction $\phi : \TT_\ep(\La)_{\geq} \to \TT_\ep(M_\ep)_{\geq}$ is an $\Abbe$-algebra isomorphism.

A pair $(M_{\ep}, \B)$ (consisting of a root of unity toric frame and an exchange matrix) 
is called a {\em root of unity quantum seed} if $(\La_{M_{\ep}}, \B)$ is $\ell$-compatible, that is 
\[
\La^\top \ol{\B} = \begin{bmatrix}
\, \ol{D} \, \\
0
\end{bmatrix},
\]
where $D$ is a diagonal matrix with positive integral diagonal entries which skew-symmetrizes the principal part of $\B$. For an integer matrix $Y$, we denote by $\ol{Y}$ the reduction 
of entries to $ \Zset/\ell$. The diagonal entries of $D$ are not required to be coprime. 

The {\em{seed mutation}} $\mu_k(M_{\ep}, \B):= (\mu_k(M_{\ep}), \mu_k(\B))$ in the direction of $k\in\ex$ is defined so that $\mu_k(\B)$ is given by \eqref{mukB} and 
\[ 
\mu_k(M_{\ep})(e_i) := \begin{cases} M_{\ep}(e_i) \quad &\text{if} \ i\neq k 
\\
M_{\ep}(-e_k + [b^k]_+)  +   M_{\ep}(-e_k -  [b^k]_-)
&\text{if} \ i= k.  
\end{cases}
\]
Here $e_1, \ldots, e_N$ denote the standard basis elements of $\Zset^N$ and for $c := \sum_i a_i e_i \in \Zset^N$,
\[
[c]_+ := \sum_{i : a_i \geq 0} a_i e_i, \quad [c]_- := \sum_{i : a_i \leq 0} a_i e_i.
\]
The $k$-th column of the matrix $\B$ is denoted by $b^k$. 

Mutation is involutive. Two seeds $(M_{\ep}, \B)$ and $(M_{\ep}'', \B'')$ are {\em{mutation-equivalent}}, denoted $(M_{\ep}, \B) \sim (M_{\ep}'', \B'')$, 
if one is obtained from the other by a sequence of mutations.

\bde{qca-root-unity} \hfill
\begin{enumerate}
\item
The {\em{root of unity quantum cluster algebra}} $\AA_{\ep}(M_\ep, \B, \inv)_{\kk}$ is the $\Abbe$-subalgebra of $\FF_\ep$ 
\[
\AA_\ep(M_\ep, \B, \inv)_{\kk}:=
\Abbe \lcor \hspace{1pt} M_\ep''(e_i), \hspace{1pt} M_\ep(e_j)^{-1} \ | \ \hspace{1pt} i \in [1,N], \hspace{1pt} j\in \inv, (M_\ep'', \B'') \sim (M_\ep, \B) \hspace{1pt}\rcor. 
\]

\item
The corresponding {\em{root of unity upper quantum cluster algebra}} is the $\Abbe$-subalgebra of $\FF_\ep$ 
\[
\UU_\ep(M_\ep, \B, \inv)_{\kk} := \hspace{-15pt} \bigcap_{ \ \ (M''_\ep, \B'') \sim (M_\ep, \B)} \hspace{-15pt} \TT_\ep(M''_\ep)_{\geq}.
\]
\end{enumerate}
\ede
It is clear that in the quantum root of unity situation, we have 
\begin{equation}
\label{base-change-ep}
\AA_\ep(M_\ep, \B, \inv)_{\kk} \cong \Abbe \otimes_{\Zset[\ep^{1/2}]} \AA_\ep(M_\ep, \B, \inv)_{\Zset[\ep^{1/2}]}
\end{equation}
but a similar fact for $\UU_\ep(M_\ep, \B, \inv)_{\kk}$ is unknown, unless we are in the case when $\kk$ is a finite extension of $\Zset$. 

{\em{We will use the notations $\AA_\ep(M_\ep, \B, \inv)$ and $\UU_\ep(M_\ep, \B, \inv)$ for brevity but $\kk$ will 
be an arbitrary integral domain of characteristic 0 and not just $\Zset$.}}

By the root of unity quantum Laurent Phenomenon \cite[Theorem 3.10]{NTY},
\[
\AA_\ep(M_\ep, \B, \inv) \subseteq \UU_\ep(M_\ep,  \B, \inv).
\]

\subsection{Canonical central subalgebras} For every toric frame $M_\ep$ of $\FF_\ep$ and $1 \leq i \leq N$, 
\[
M(e_i)^\ell \in \ZZ(\FF_\ep). 
\] 
By \cite[Proposition 4.4]{NTY}, if the following condition holds
\medskip

{\bf{(Coprime)}} $\ell \in \Zset_+$ is odd and coprime to the diagonal entries of the skew-symmetrizing matrix $D$, 
\medskip
\\
then for all $k \in \ex$,
\begin{equation}
\label{l-mutation}
M_\ep(e_k)^\ell \left( \mu_k M_\ep (e_k) \right)^{\ell} = \prod_{b_{ik}>0} (M_\ep(e_i)^\ell)^{b_{ik}} + \prod_{b_{ik}<0} (M_\ep(e_i)^\ell)^{ - b_{ik}}.
\end{equation}
By Theorem 4.6 and Corollary 4.7 in \cite{NTY}, the central  $\Abbe$-subalgebra of $\AA_\ep(M_\ep, \B, \inv)$ 
\[
\CC_\ep(M_\ep, \B, \inv):=
\Abbe \lcor \hspace{1pt} M_\ep''(e_i)^\ell, \hspace{1pt} M_\ep(e_j)^{-\ell} \ | \ \hspace{1pt} i \in [1,N], \hspace{1pt} j\in \inv, (M_\ep'', \B'') \sim (M_\ep, \B) \hspace{1pt}\rcor  
\]
is isomorphic to a base change of the underlying cluster algebra:
\[
\CC_\ep(M_\ep, \B, \inv)\cong \Abbe \otimes_{\kk} \AA(\B, \inv).
\]
%%%%%%
\subsection{Exchange graphs}
The {\em{exchange graphs}} of an upper cluster algebra $\UU(\B, \inv)$ and a root of unity upper quantum cluster algebra $\UU_\ep(M_\ep, \B, \inv)$ 
are the labelled graphs with vertices corresponding to the seeds that are mutation-equivalent to $(\wt{\mathbf{x}}, \B)$ and $(M_\ep, \B)$, respectively,
and edges given by seed mutations and labelled by the corresponding mutation numbers.
The exchange graphs are independent on the choice to work with cluster algebras vs their upper counterparts,
the choice of base ring $\kk$ by \eqref{base-change} and \eqref{base-change-ep}, 
and the choice of the inverted set of frozen indices $\inv$. 
Those graphs will be denoted by $E(\B)$ and $E_\ep(\La_{M_\ep},\B)$, respectively.

By \cite[Theorem 4.8]{NTY}, if condition (Coprime) is satisfied, then there is a unique isomorphism of labelled graphs 
\begin{equation}
\label{graph-iso}
E_\ep(M_\ep, \B) \cong E(\B)
\end{equation}
sending the vertex corresponding to the seed $(M_\ep, \B)$ to the vertex corresponding to the seed  $(\wt{\mathbf{x}}, \B)$.
%%%%%%%%%%%%%
\sectionnew{Root of unity upper quantum cluster algebras}
%%%%%%%%%%
This section contains material on root of unity upper quantum cluster algebras, 
partial intersections of mixed quantum tori, their special and full centers.

In the rest of the paper we retain the notation from Sect. \ref{2.3}, and $\kk$ will denote an
arbitrary integral domain of characteristic 0.
\subsection{The algebras $\UU_\ep(M_\ep, \B, \inv)$ and their centers}
For a root of unity toric frame $M_\ep$ of an $\Abbe$-division algebra $\FF_\ep$, the subfields of 
$\FF_\ep$ generated by $\{M_\ep(e_i)^\ell \mid 1 \leq i \leq N\}$ and $\Abbe \cup \{M_\ep(e_i)^\ell \mid 1 \leq i \leq N\}$ 
are purely transcendental extensions of $F$ and $F(\ep)$, respectively, of degree $N$. (Recall that $F$ denotes the fraction field of the integral domain $\kk$.)
They contain the mixed  polynomial/Laurent polynomial rings
\begin{align}
\TT(M_\ep)^\ell_{\geq} &:= \kk[M_\ep(e_i)^{\pm \ell}, M_\ep(e_j)^\ell; i \in \ex \sqcup \inv, j \notin \ex \sqcup \inv], 
\label{TlZ}
\\
\TT_\ep(M_\ep)^\ell_{\geq} &:= \Abbe[M_\ep(e_i)^{\pm \ell}, M_\ep(e_j)^\ell; i \in \ex \sqcup \inv, j \notin \ex \sqcup \inv].
\label{Tl}
\end{align}
Clearly, $\TT_\ep(M_\ep)^\ell_{\geq} \cong \TT(M_\ep)^\ell_{\geq} \otimes_{\kk} \Abbe$.
Denote the following central subring and subalgebra of $\FF_\ep$:
\begin{align*}
\CC\UU(M_\ep, \B, \inv) &:= \hspace{-15pt} \bigcap_{ \ \ (M''_\ep, \B'') \sim (M_\ep, \B)} \hspace{-15pt} \TT(M''_\ep)^\ell_{\geq}, \\
\CC\UU_\ep(M_\ep, \B, \inv) &:= \hspace{-15pt} \bigcap_{ \ \ (M''_\ep, \B'') \sim (M_\ep, \B)} \hspace{-15pt} \TT_\ep(M''_\ep)^\ell_{\geq}.
\end{align*}
As in the previous section, for brevity, we will not display the base ring $\kk$ 
in the notations $\CC\UU(M_\ep, \B, \inv)$ and  $\CC\UU_\ep(M_\ep, \B, \inv)$.

%If condition (Coprime) holds, then one easily sees that \eqref{l-mutation} implies that
%\[
%\CC\UU_\ep(M_\ep, \B, \inv) \cong \CC\UU(M_\ep, \B, \inv) \otimes_\Zset \Abbe.
%\]
\bpr{cent-U} Assume that $\ell \in \Zset_+$ and $\B \in M_{N \times \ex}(\Zset)$ is an exchange matrix such that 
$\ell$ satisfies condition (Coprime). Then 
\[
\CC\UU(M_\ep, \B, \inv) \cong \UU(\B, \inv)_{\kk} \quad \mbox{and} \quad \CC\UU_\ep(M_\ep, \B, \inv) \cong \UU(\B, \inv)_{\Abe}.
\]
\epr
\begin{proof}
It follows from \eqref{l-mutation} that for all seeds $(M''_\ep, \B'') \sim (M_\ep, \B)$, the mixed polynomial/Laurent polynomial ring
$\TT(M''_\ep)^\ell_{\geq}$ lies inside the fraction field of $\TT(M_\ep)^\ell_{\geq}$. Furthermore, \eqref{l-mutation} also implies 
that the generators of $\TT(M''_\ep)^\ell_{\geq}$ obey the classical mutation rule \eqref{classical-mutation}. The statement of the 
proposition now follows from the isomorphism \eqref{graph-iso} of the exchange graphs of the classical cluster algebra 
and its root of unity quantum counterpart. 
\end{proof}
%%%%%%%%%%%%
\subsection{Partial intersections}
For the purposes of flexibility of application to representation theory, we consider partial intersections of mixed quantum tori that 
generalize root of unity upper quantum cluster algebras. For a subset $\Theta$ of seeds, denote the $\Abbe$-algebra
\[
\UU_\ep(M_\ep, \B, \inv, \Theta) := \hspace{-15pt} \bigcap_{ \ \ (M''_\ep, \B'') \in \Theta} \hspace{-15pt} \TT_\ep(M''_\ep)_{\geq}.
\]
The following subring and subalgebra of $\UU_\ep(M_\ep, \B, \inv, \Theta)$
\begin{align*}
\CC\UU(M_\ep, \B, \inv, \Theta) &:= \hspace{-15pt} \bigcap_{ \ \ (M''_\ep, \B'') \in \Theta} \hspace{-15pt} \TT(M''_\ep)^\ell_{\geq}, \\
\CC\UU_\ep(M_\ep, \B, \inv, \Theta) &:= \hspace{-15pt} \bigcap_{ \ \ (M''_\ep, \B'') \in \Theta} \hspace{-15pt} \TT_\ep(M''_\ep)^\ell_{\geq}
\end{align*}
lie in its center and $\CC\UU_\ep(M_\ep, \B, \inv, \Theta) \cong \CC\UU(M_\ep, \B, \inv, \Theta) \otimes_{\kk} \Abbe$.

Analogously to \prref{cent-U} one proves  the following:
\bpr{cent-U-Theta} If $\ell$ satisfies condition (Coprime), then 
\[
\CC\UU(M_\ep, \B, \inv, \Theta) \cong \UU(\B, \inv, \Theta)_{\kk} := \bigcap_{(\wt{\mathbf{x}}'', \B'') \in \Theta} \kk(x''_i)^{\pm 1}, x''_j; i \in \ex \sqcup \inv, j \notin \ex \sqcup \inv],
\]
where in the intersection we use the isomorphism \eqref{graph-iso} to identify $\Theta$ with a subset of the set of vertices of the exchange graph of 
$\UU(\B, \inv)$, and 
\[
\CC\UU_\ep(M_\ep, \B, \inv, \Theta) \cong \UU(\B, \inv, \Theta)_{\Abbe}.
\]
\epr
\subsection{Full centers} 
\ble{centU}
For any subset $\Theta$ of the set of vertices of the exchange graph $E_\ep(M_\ep, \B)$ of the root of unity upper quantum cluster algebra 
$\UU_\ep(M_\ep, \B, \inv)$ (which is not necessarily connected) and any root of unity 
$\ep$ (without restrictions on its order), the center of $\UU_\ep(M_\ep, \B, \inv, \Theta)$ is given by 
\begin{equation}
\label{full-cent}
\ZZ(\UU_\ep(M_\ep, \B, \inv, \Theta)) = \UU_\ep(M_\ep, \B, \inv, \Theta) \cap \ZZ(\TT_\ep(M'_\ep))
\end{equation}
for any seed $(M'_\ep, \B') \in \Theta$. 

In particular, 
\[
\ZZ(\UU_\ep(M_\ep, \B, \inv)) = \UU_\ep(M_\ep, \B, \inv) \cap \ZZ(\TT_\ep(M'_\ep))
\]
for any seed $(M'_\ep, \B')$ of $\UU_\ep(M_\ep, \B, \inv)$.
\ele
\begin{proof}
Eq. \eqref{full-cent} follows at once from the fact that $\TT_\ep(M'_\ep)$ is a central localization of the algebra $\UU_\ep(M_\ep, \B, \inv, \Theta)$:
\begin{equation}
\label{TU}
\TT_\ep(M'_\ep) \cong \UU_\ep(M_\ep, \B, \inv, \Theta)[ M'_\ep(e_i)^{-\ell}, 1 \leq i \leq N].
\end{equation}
\end{proof}

The center $\ZZ(\TT_\ep(M'_\ep))$ of the quantum torus $\TT_\ep(M'_\ep)$ is explicitly described in Sect. \ref{6.1}. In Sect. \ref{4.3} we present 
similar descriptions of the central subalgebras $\CC\UU(M_\ep, \B, \inv, \Theta)$ and $\CC\UU_\ep(M_\ep, \B, \inv, \Theta)$ as intersections
of the form \eqref{full-cent}.
%%%%%%%%%%%%%
\sectionnew{A Cayley--Hamilton structure on $(\UU_\ep(M_\ep, \B, \inv, \Theta), \CC\UU_\ep(M_\ep, \B, \inv, \Theta))$}
%%%%
In this section we construct Cayley--Hamilton structures on root of unity upper quantum cluster algebras with respect to their special centers,
proving parts (1)-(4) of Theorem B from the introduction. 
\subsection{Mixed quantum tori}
\label{4.1}
If $R$ is a $k$-algebra over a commutative ring $k$ which is free over a central subalgebra $C$ of finite rank $r$, then the left action of $R$ on itself induces the algebra homomorphism 
\[
R \to \End_{C}(R) \cong M_r(C). 
\]
The {\em{regular trace}} of $R$ with respect to $C$ is the composition of this map with the matrix trace on $M_r(C)$: 
\[
\tr_{\reg}^{R,C}: R \to \End_C(R) \cong M_r(C) \stackrel{\Tr}{\to} C.
\]
If $i$ is not a zero divisor of $R$ for $1 \leq i \leq r$, then
\medskip

(*) the triple $(R, C, \tr_{reg}^{R,C})$ is a Cayley--Hamilton algebra of degree $r$.
\medskip
\\
This follows from the fact that $(M_r(C), C, \Tr)$ is a Cayley--Hamilton algebra of degree $r$.

Recall the notation \eqref{Tl}.

\ble{T-Tl} The following hold for a root of unity toric frame $M_\ep$ of an $\Abbe$-division algebra $\FF_\ep$:
\begin{enumerate}
\item $\TT_\ep( M_\ep) _{\geq}$ is a free $\TT_\ep(M_\ep)_{\geq}^\ell$-module of rank $\ell^N$ with basis
\[
\{ \hspace{1pt} M_\ep (f) \hspace{1pt} | \hspace{1pt} f=(f_1, f_2, \dots, f_N)\in \Zset_{\geq}^N, 0\le f_i< \ell \hspace{1pt} \}.
\]
\item The regular trace for the pair $(\TT_\ep( M_\ep) _{\geq}, \TT_\ep(M_\ep)_{\geq}^\ell)$, to be denoted by $\tr_{\reg}^{\TT_\ep(M_\ep)_{\geq}}$, 
is given by 
\[
\tr_{\reg}^{\TT_\ep(M_\ep)_{\geq}}(M_\ep(f)) = 
\begin{cases}
\ell^N\cdot M_\ep (f), & \text{if } f \in (\ell\Zset)^N\\
0,                  &      \text{if } f\not\in (\ell\Zset)^N.
\end{cases}
\]
\item The triple $(\TT_\ep( M_\ep) _{\geq}, \TT_\ep(M_\ep)_{\geq}^\ell, \tr_{\reg}^{\TT_\ep(M_\ep)_{\geq}})$ is a Cayley--Hamilton algebra
of degree $\ell^N$.
\end{enumerate}
\ele
\begin{proof} Part (1) follows from the product formula \eqref{prod-quant-tor}. The second part is straightforward. The third is a special case of statement (*) above.
\end{proof}
%%%%%%%%%%
\subsection{Structure of neighboring intersections}
Consider a root of unity quantum seed $\ (M_\ep, \B)$ and a 1-step mutation $\mu_k(M_{\ep}, \B)= (M'_{\ep}, \B')$. 
In this subsection we analyze the structure of the intersection
\[
R:=\TT_\ep(M_\ep)_{\geq} \cap \TT_\ep(M'_\ep)_{\geq}. 
\]

Denote 
\[
x_k:=M_\ep(e_k) \quad \mbox{and} \quad y_k:=M'_\ep(e_k).
\]
Define the mixed quantum torus
\[
\TT_\ep (M_\ep)^\circ_{\geq}:= \Abbe\hspace{-0.1cm}-\hspace{-0.08cm}\Span \{M_\ep(f) \mid f=(f_1, f_2, \dots, f_N) \in \Zset^N, f_k=0, f_j\ge 0 \text{ for }  j\not\in \ex \sqcup \inv \}. 
\]
sitting inside $R$. It has the direct sum decomposition 
\[
\TT_\ep (M_\ep)^\circ_{\geq} = \TT_\ep (M_\ep)^{\circ, \cc}_{\geq} \oplus \TT_\ep (M_\ep)^{\circ, \nc}_{\geq},
\]
where
\begin{align*}
\TT_\ep (M_\ep)^{\circ, \cc}_{\geq}&:= \Abbe\hspace{-0.1cm}-\hspace{-0.08cm}\Span \{M_\ep(f) \mid f=(f_1, f_2, \dots, f_N) \in (\ell \Zset)^N, 
\\
& \hspace{5cm} f_k=0, f_j\ge 0 \text{ for }  j\not\in \ex \sqcup \inv \}, 
\\
\TT_\ep (M_\ep)^{\circ, \nc}_{\geq}&:= \Abbe\hspace{-0.1cm}-\hspace{-0.08cm}\Span \{M_\ep(f) \mid f=(f_1, f_2, \dots, f_N) \in \Zset^N \backslash (\ell \Zset)^N, 
\\
& \hspace{5cm} f_k=0, f_j\ge 0 \text{ for }  j\not\in \ex \sqcup \inv \}.
\end{align*}
The algebra $\TT_\ep (M_\ep)^{\circ, \cc}_{\geq}$ is inside $\ZZ(R)$ and $\TT_\ep (M_\ep)^{\circ, \nc}_{\geq}$ is a $\TT_\ep (M_\ep)^{\circ, \cc}_{\geq}$-module
under the left (and right) action.

The mixed quantum torus $\TT_\ep(M_\ep)_{\geq}$ is decomposed as
\begin{equation}
\label{decomp1}
\TT_\ep(M_\ep)_{\geq} = \Big( \bigoplus_{n \in \Zset, \ell \nmid n} x_k^n \TT_\ep(M_\ep)_{\geq}^{\circ} \Big) \oplus 
\Big( \bigoplus_{n \in \Zset} x_k^{\ell n} \TT_\ep(M_\ep)_{\geq}^{\circ, \nc} \Big) \oplus
\Big( \bigoplus_{n \in \Zset} x_k^{\ell n} \TT_\ep(M_\ep)_{\geq}^{\circ, \cc} \Big).
\end{equation}
The third term is $\TT_\ep(M_\ep)_\geq^\ell$:
\begin{equation}
\label{Tl-pf1}
\TT_\ep(M_\ep)^\ell_\geq = \bigoplus_{n \in \Zset} x_k^{\ell n} \TT_\ep (M_\ep)^{\circ, \cc}_{\geq}.
\end{equation}
By \leref{T-Tl}(2), 
\medskip

(**) the trace function $\tr_{\reg}^{\TT_\ep(M_\ep)_{\geq}}$ vanishes on the first two summands and equals $\ell^N \cdot \Id$ on the third one. 
\medskip
\\
Likewise the mixed quantum torus $\TT_\ep(M'_\ep)_{\geq}$ is decomposed as
\begin{equation}
\label{decomp2}
\TT_\ep(M'_\ep)_{\geq} = \Big( \bigoplus_{n \in \Zset, \ell \nmid n} y_k^n \TT_\ep(M_\ep)_{\geq}^{\circ} \Big) \oplus 
\Big( \bigoplus_{\ell n \in \Zset} y_k^{\ell n} \TT_\ep(M_\ep)_{\geq}^{\circ, \nc} \Big) \oplus
\Big( \bigoplus_{\ell n \in \Zset} y_k^{\ell n} \TT_\ep(M_\ep)_{\geq}^{\circ, \cc} \Big),
\end{equation}
where the third term is $\TT_\ep(M'_\ep)_\geq^\ell$:
\begin{equation}
\label{Tl-pf2}
\TT_\ep(M'_\ep)^\ell_\geq = \bigoplus_{n \in \Zset} y_k^{\ell n} \TT_\ep (M_\ep)^{\circ, \cc}_{\geq}.
\end{equation}
By \leref{T-Tl}(2), 
\medskip

(***) the trace function $\tr_{\reg}^{\TT_\ep(M'_\ep)_{\geq}}$ vanishes on the first two summands and equals $\ell^N \cdot \Id$ on the third one. 

\bth{inters} Let $\ell$ satisfies condition (Coprime).
Assume that $\mu_k(M_{\ep}, \B)= (M'_{\ep}, \B')$. 
For all $n \in \Zset$, we have:
\begin{enumerate}
\item
$x_k^n \TT_\ep(M_\ep)_{\geq}^\circ \bigcap \TT_\ep(M'_\ep)_{\geq}= \TT_\ep(M_\ep)_{\geq}\bigcap y_k^{-n} \TT_\ep(M_\ep)_{\geq}^{\circ}= 
x_k^{n}\TT_\ep(M_\ep)_{\geq}^{\circ} \bigcap y_k^{-n} \TT_\ep(M_\ep)_{\geq}^{\circ}$. \\
\item
$x_k^{\ell n} \TT_\ep(M_\ep)_{\geq}^{\circ, \cc} \bigcap \TT_\ep(M'_\ep)_{\geq}= 
\TT_\ep(M_\ep)_{\geq}\bigcap y_k^{- \ell n} \TT_\ep(M_\ep)_{\geq}^{\circ, \cc} =
x_k^{\ell n} \TT_\ep(M_\ep)_{\geq}^{\circ, \cc} \bigcap y_k^{- \ell n} \TT_\ep(M_\ep)_{\geq}^{\circ, \cc};$ \\
\item
$\TT_\ep(M_\ep)_{\geq}^\ell \bigcap \TT_\ep(M'_\ep)_{\geq}= 
\TT_\ep(M_\ep)_{\geq}\bigcap \TT_\ep(M'_\ep)_{\geq}^\ell =
\TT_\ep(M_\ep)_{\geq}^\ell \bigcap \TT_\ep(M'_\ep)_{\geq}^\ell;$ \\
\item
$x_k^{\ell n} \TT_\ep(M_\ep)_{\geq}^{\circ, \nc} \bigcap \TT_\ep(M'_\ep)_{\geq}= 
\TT_\ep(M_\ep)_{\geq}\bigcap y_k^{- \ell n} \TT_\ep(M_\ep)_{\geq}^{\circ, \nc} =
x_k^{\ell n} \TT_\ep(M_\ep)_{\geq}^{\circ, \nc} \bigcap y_k^{- \ell n} \TT_\ep(M_\ep)_{\geq}^{\circ, \nc};$ \\
\item The direct sum decompositions \eqref{decomp1} and \eqref{decomp2} restrict to direct sum decompositions of $R$. 
\item
The regular trace maps $\tr_{\reg}^{\TT_\ep(M_\ep)_{\geq}}$ and $\tr_{\reg}^{\TT_\ep(M'_\ep)_{\geq}}$ 
coincide on $R$ and map it to \\ $\TT_\ep(M_\ep)^\ell_{\geq}\bigcap\TT_\ep(M'_\ep)^\ell_{\geq}$.
\end{enumerate}
\eth
\begin{proof} We use the idea of \cite[Proposition 3.5]{GLS-spec} to deal with the 
elements of $R$. Define 
\[
Q_n = \ep^{n \La(e_k, [b^k]_+)/2} M_\ep([b^k]_+) + \ep^{- n \La(e_k, [b^k]_-)/2}  M_\ep(- [b^k]_-)  \in \TT_\ep(M_\ep)^\circ_{\geq}
\]
for $n \in \Zset$. Then 
\[
Q_1 = x_k y_k\quad \mbox{and} \quad Q_n x_k = x_k Q_{n-2}, \; \forall n \in \Zset.
\]
All elements $r \in R$ are of the form
\begin{equation}
\label{r}
r = \sum_{n \in \Zset} x_k^n a_n  = \sum_{n \in \Zset} (y_k)^n c_n,   
\end{equation}
where both sums are finite and $a_n, c_n \in \TT_\ep(M_\ep)^\circ_{\geq}$, $\forall n \in \Zset$. 
Since the mixed quantum torus $\TT_\ep(M_\ep)_\geq$ has the basis $\{M(f) \mid f \in \Zset^N, f_j \geq 0, \forall j \notin \ex \sqcup \inv\}$,
\begin{align}
& a_0 =c_0, &&
\nn
\\
&a_n = Q_{- 2 n -1} \ldots Q_3 Q_1 c_{-n}, \quad &&\forall n < 0, 
\nn
\\
&Q_{-1} Q_{-3} \ldots Q_{-2n+1} a_n = c_{-n}, \quad &&\forall n > 0.
\nn
\end{align} 

(1) If $r \in x_k^n \TT_\ep(M_\ep)_{\geq}^\circ \bigcap \TT_\ep(M'_\ep)_{\geq}$, then $a_m = 0$ for all $m \neq n$. Hence $c_m=0$ for
$m \neq -n$, and thus, $r \in x_k^{n}\TT_\ep(M_\ep)_{\geq}^{\circ} \bigcap y_k^{-n} \TT_\ep(M_\ep)_{\geq}^{\circ}$. Since
\[
x_k^{n}\TT_\ep(M_\ep)_{\geq}^{\circ} \bigcap y_k^{-n} \TT_\ep(M_\ep)_{\geq}^{\circ} \subseteq x_k^n \TT_\ep(M_\ep)_{\geq}^\circ \bigcap \TT_\ep(M'_\ep)_{\geq}, 
\]
we get that 
\[
x_k^n \TT_\ep(M_\ep)_{\geq}^\circ \bigcap \TT_\ep(M'_\ep)_{\geq} = x_k^{n}\TT_\ep(M_\ep)_{\geq}^{\circ} \bigcap y_k^{-n} \TT_\ep(M_\ep)_{\geq}^{\circ}.
\]
By interchanging the roles of $M_\ep$ and $M'_\ep$, we obtain
\[
\TT_\ep(M_\ep)_{\geq}\bigcap y_k^{-n} \TT_\ep(M_\ep)_{\geq}^{\circ} =
x_k^{n}\TT_\ep(M_\ep)_{\geq}^{\circ} \bigcap y_k^{-n} \TT_\ep(M_\ep)_{\geq}^{\circ}.
\]

For parts (2) and (3), denote 
\[
Q:=\prod_{b_{ik}>0} (M_\ep(e_i)^\ell)^{b_{ik}} + \prod_{b_{ik}<0} (M_\ep(e_i)^\ell)^{ - b_{ik}} \in \TT_\ep(M_\ep)_{\geq}^{\circ, \cc}
\]
and note that \eqref{l-mutation} implies that
\[
x_k^{\ell n} y_k^{\ell n} = Q^n, \quad \forall n \geq 0.
\]
Therefore the coefficients $a_{\ell n}$ and $c_{\ell n}$ in \eqref{r} satisfy
\[
a_{-\ell n} = Q^n c_{\ell n}, \quad Q^n a_{\ell n} = c_{- \ell n} \quad \forall n \geq 0.
\]

(2) If $r \in x_k^{\ell n} \TT_\ep(M_\ep)_{\geq}^{\circ, \cc} \bigcap \TT_\ep(M'_\ep)_{\geq}$, 
then $a_{\ell n} \in  \TT_\ep(M_\ep)_{\geq}^{\circ, \cc}$ and by part (1), $c_m = 0$ for 
$m \neq - \ell n$. 

Case 1: If $n \geq 0$, then 
\[
c_{- \ell n} = Q^n a_{\ell n} \in \TT_\ep(M_\ep)_{\geq}^{\circ, \cc}
\]
because $Q^n, a_{\ell n} \in \TT_\ep(M_\ep)_{\geq}^{\circ, \cc}$. Therefore, 
$r \in x_k^{\ell n} \TT_\ep(M_\ep)_{\geq}^{\circ, \cc} \bigcap y_k^{- \ell n} \TT_\ep(M_\ep)_{\geq}^{\circ, \cc}$.

Case 2: If $n < 0$, then 
\[
c_{- \ell n} = Q^n a_{\ell n} \in \TT_\ep(M_\ep)_{\geq}^{\circ, \cc}[Q^{-1}] \cap \TT_\ep(M_\ep)_{\geq}^{\circ} = \TT_\ep(M_\ep)_{\geq}^{\circ, \cc}
\]
and once again $r \in x_k^{\ell n} \TT_\ep(M_\ep)_{\geq}^{\circ, \cc} \bigcap y_k^{- \ell n} \TT_\ep(M_\ep)_{\geq}^{\circ, \cc}$.

Combining the two cases gives
\[
x_k^{\ell n} \TT_\ep(M_\ep)_{\geq}^{\circ, \cc} \bigcap \TT_\ep(M'_\ep)_{\geq} \subseteq 
x_k^{\ell n} \TT_\ep(M_\ep)_{\geq}^{\circ, \cc} \bigcap y_k^{- \ell n} \TT_\ep(M_\ep)_{\geq}^{\circ, \cc}
\]
and the opposite inclusion is obvious. Hence, 
\[
x_k^{\ell n} \TT_\ep(M_\ep)_{\geq}^{\circ, \cc} \bigcap \TT_\ep(M'_\ep)_{\geq} =
x_k^{\ell n} \TT_\ep(M_\ep)_{\geq}^{\circ, \cc} \bigcap y_k^{- \ell n} \TT_\ep(M_\ep)_{\geq}^{\circ, \cc}.
\]
The equality
\[
\TT_\ep(M_\ep)_{\geq}\bigcap y_k^{- \ell n} \TT_\ep(M_\ep)_{\geq}^{\circ, \cc} =
x_k^{\ell n} \TT_\ep(M_\ep)_{\geq}^{\circ, \cc} \bigcap y_k^{- \ell n} \TT_\ep(M_\ep)_{\geq}^{\circ, \cc}
\]
is proved by interchanging the roles of $M_\ep$ and $M'_\ep$. 

Part (3) follows from part (1) and Eqs. \eqref{Tl-pf1} and \eqref{Tl-pf2}.
Part (4) is proved similarly to part (2). Part (5) is proved by combining the arguments in the proofs of parts (1), (2) and (4).

(6) For an element $r \in R$, denote by $r_1, r_2, r_3$ its components in the direct sum decomposition \eqref{decomp1} 
and by $s_1, s_2, s_3$ its components in the direct sum decomposition \eqref{decomp2}. By facts (**) and (***) above, 
\begin{equation}
\label{tr-reg-12}
\tr_{\reg}^{\TT_\ep(M_\ep)_{\geq}}(r) = \ell^N r_3 \quad \mbox{and} \quad \tr_{\reg}^{\TT_\ep(M'_\ep)_{\geq}}(r) = \ell^N s_3.
\end{equation}
Parts (1), (2), (4) and (5) of the theorem imply that 
\begin{align*}
&r_1, s_1 \in \Big( \bigoplus_{n \in \Zset, \ell \nmid n} x_k^n \TT_\ep(M_\ep)_{\geq}^{\circ} \bigcap y_k^{- n} \TT_\ep(M_\ep)_{\geq}^{\circ} \Big),
\\
& r_2, s_2 \in \Big( \bigoplus_{n \in \Zset} x_k^{\ell n} \TT_\ep(M_\ep)_{\geq}^{\circ, \nc} \bigcap  y_k^{-\ell n} \TT_\ep(M_\ep)_{\geq}^{\circ, \nc}\Big),
\\
& r_3, s_3 \in \Big( \bigoplus_{n \in \Zset} x_k^{\ell n} \TT_\ep(M_\ep)_{\geq}^{\circ, \cc} \bigcap  y_k^{-\ell n} \TT_\ep(M_\ep)_{\geq}^{\circ, \cc}\Big).
\end{align*}
Therefore $r_1 = s_1$, $r_2 = s_2$, $r_3 = s_3$, and by \eqref{tr-reg-12}, 
$\tr_{\reg}^{\TT_\ep(M_\ep)_{\geq}}(r) = \tr_{\reg}^{\TT_\ep(M'_\ep)_{\geq}}(r)$. 
\end{proof}
\subsection{The central subrings $\CC\UU(M_\ep, \B, \inv, \Theta)$ and $\CC\UU_\ep(M_\ep, \B, \inv, \Theta)$}
\label{4.3}
\thref{inters} leads to the following description of the two canonical central subring and subalgebra of $\UU_\ep(M_\ep, \B, \inv, \Theta)$
in a manner that is similar to the description of the full center of $\UU_\ep(M_\ep, \B, \inv, \Theta)$ from \leref{centU}.
\bpr{CU} Assume that $\ell$ satisfies condition (Coprime). For every connected subset 
$\Theta$ of the exchange graph $E_\ep(\La_{M_\ep},\B)$ of $\UU(M_\ep, \B, \inv)$ and any seed
$(M'_\ep, \B')$ in $\Theta$, we have
\begin{equation}
\label{CU1}
\CC\UU(M_\ep, \B, \inv, \Theta) = \UU_\ep(M_\ep, \B, \inv, \Theta) \cap \TT(M'_\ep)^\ell_{\geq}
\end{equation}
and
\begin{equation}
\label{CU2}
\CC\UU_\ep(M_\ep, \B, \inv, \Theta) = \UU_\ep(M_\ep, \B, \inv, \Theta) \cap \TT_\ep(M'_\ep)^\ell_{\geq},
\end{equation}
recall \eqref{TlZ} and \eqref{Tl}.

\epr
\begin{proof} First we show \eqref{CU2}. For two seeds $(M'_\ep, \B'), (M''_\ep, \B'') \in \Theta$
denote by 
\[
d \big( (M'_\ep, \B'), (M''_\ep, \B'') \big)
\]
the distance between them in the full subgraph of $E_\ep(\La_{M'_\ep},\B)$ with vertex set $\Theta$. For $k \geq 0$, denote the 
subalgebras 
\begin{align*}
C_k := &\big( \cap \{ \TT_\ep(M''_\ep)^\ell_{\geq} \mid (M''_\ep, \B'') \in \Theta, d( (M'_\ep, \B'), (M''_\ep, \B'')) \leq k \} \big) \bigcap
\\
&\big( \cap \{ \TT_\ep(M''_\ep)_{\geq} \mid (M''_\ep, \B'') \in \Theta, d( (M'_\ep, \B'), (M''_\ep, \B'')) > k \} \big).
\end{align*}
It is clear that
\[
\UU_\ep(M_\ep, \B, \inv, \Theta) \cap \TT_\ep(M'_\ep)^\ell_{\geq} = C_ 0 \supseteq C_1 \supseteq \ldots
\]
and that
\[
\bigcap_{k \geq 0} C_k = \CC\UU_\ep(M_\ep, \B, \inv, \Theta).
\]
By recursively applying \thref{inters}(3), we obtain that
\[
C_0 = C_1 = \ldots
\]
This proves \eqref{CU2}.

Analogously to \thref{inters}(3), one proves that, if $\mu_k(M_{\ep}, \B)= (M'_{\ep}, \B')$, then
\[
\TT(M_\ep)_{\geq}^\ell \bigcap \TT_\ep(M'_\ep)_{\geq}= 
\TT_\ep(M_\ep)_{\geq}\bigcap \TT(M'_\ep)_{\geq}^\ell =
\TT(M_\ep)_{\geq}^\ell \bigcap \TT(M'_\ep)_{\geq}^\ell.
\]
Eq. \eqref{CU1} is deduced from this property by an analogous argument to the one used for \eqref{CU2}.
\end{proof}
\subsection{Construction of Cayley--Hamilton structures} 
\bth{CH-U-CU}
Let $\B$ be an exchange matrix and $\ep^{1/2}$ be a primitive $\ell$-th root of unity for an integer $\ell$ that satisfies condition (Coprime). Let $\Theta$ be a connected set of vertices of the exchange graph 
$E_\ep(\La_{M_\ep},\B)$ of $\UU_\ep(M_\ep, \B, \inv)$ (recall the isomorphism \eqref{graph-iso} to the 
exchange graph $E(\B)$ of $\UU(\B, \inv)$). The following hold:
\begin{enumerate}
\item
For every pair of seeds $(M'_\ep, \B'), (M''_\ep, \B'') \in \Theta$, 
\begin{equation}
\label{t-eq-t}
\tr_{\reg}^{\TT_\ep(M'_\ep)_{\geq}}\big{|}_{\UU_\ep(M_\ep, \B, \inv, \Theta)} = \tr_{\reg}^{\TT_\ep(M''_\ep)_{\geq}}\big{|}_{\UU_\ep(M_\ep, \B, \inv, \Theta)}.
\end{equation}
We denote by $\tr_{\reg}$ this restriction map coming from an arbitrary seed in $\Theta$. 
\item $\tr_{\reg} ( \UU_\ep(M_\ep, \B, \inv, \Theta)) \subseteq \CC\UU_\ep(M_\ep, \B, \inv, \Theta)$.
\item The triple 
\[
(\UU_\ep(M_\ep, \B, \inv, \Theta), \CC\UU_\ep(M_\ep, \B, \inv, \Theta), \tr_{\reg}) 
\]
is a Cayley--Hamilton algebra of degree equal to $\ell^N$.
\item Assume that the base ring $\kk$ is a field extension of the cyclotomic field $\Qset(\ep^{1/2})$.
Then $\UU_\ep(M_\ep, \B, \inv, \Theta)$ is a finitely generated $\kk$-algebra if and only 
$\CC \UU_\ep(M_\ep, \B, \inv, \Theta)$ is a finitely generated $\kk$-algebra and 
$\UU_\ep(M_\ep, \B, \inv, \Theta)$ is a finitely generated module over $\CC \UU_\ep(M_\ep, \B, \inv, \Theta)$.
\end{enumerate}
\eth
\begin{proof}
(1) Assume that $(M'_\ep, \B')$ and $(M''_\ep, \B'') := \mu_k (M'_\ep, \B')$ are two adjacent seeds in $\Theta$. Since, 
\[
\UU_\ep(M_\ep, \B, \inv, \Theta) \subseteq \TT_\ep(M'_\ep)_\geq \cap  \TT_\ep(M''_\ep)_\geq,
\]
\thref{inters}(6) implies that \eqref{t-eq-t} holds for pairs of adjacent seeds in $\Theta$. Because $\Theta$ is a connected subset of vertices of the exchange graph 
$E_\ep(\La_{M_\ep},\B)$,  \eqref{t-eq-t} holds for pairs of seeds in $\Theta$.

For part (2) we have: 
\begin{align*}
\tr_{\reg} ( \UU_\ep(M_\ep, \B, \inv, \Theta)) &\subseteq 
\bigcap_{(M'_\ep, \B') \in \Theta} \Im \tr_{\reg}^{\TT_\ep(M'_\ep)_{\geq}} = \bigcap_{(M'_\ep, \B') \in \Theta} \TT_\ep(M'_\ep)_{\geq}^\ell 
\\
&= \CC\UU_\ep(M_\ep, \B, \inv, \Theta).
\end{align*}

Part (3) follows from Lemmas \ref{lrestrict-CH} and \ref{lT-Tl}(3) and part (2) of the theorem.

(4) Assume first that $\UU_\ep(M_\ep, \B, \inv, \Theta)$ is a finitely generated $\kk$-algebra. This assumption and the 
Cayley--Hamilton algebra structure from part (3) make possible the application of \cite[Theorem 4.5]{DP}, 
which implies that $\CC \UU_\ep(M_\ep, \B, \inv, \Theta)$ is a finitely generated $\kk$-algebra and 
$\UU_\ep(M_\ep, \B, \inv, \Theta)$ is a finitely generated module over $\CC \UU_\ep(M_\ep, \B, \inv, \Theta)$.

In the opposite direction, if $\CC \UU_\ep(M_\ep, \B, \inv, \Theta)$ is a finitely generated $\kk$-algebra and 
$\UU_\ep(M_\ep, \B, \inv, \Theta)$ is a finitely generated module over $\CC \UU_\ep(M_\ep, \B, \inv, \Theta)$, 
then,  as a $\kk$-algebra, $\UU_\ep(M_\ep, \B, \inv, \Theta)$ is generated by the collection of generators of the $\kk$-algebra 
$\CC \UU_\ep(M_\ep, \B, \inv, \Theta)$ together with the collection of generators 
of $\UU_\ep(M_\ep, \B, \inv, \Theta)$ as a $\CC \UU_\ep(M_\ep, \B, \inv, \Theta)$-module.
\end{proof}
\bre{reg-tr-U} Although in \thref{CH-U-CU} we denote the trace function 
\[
\UU_\ep(M_\ep, \B, \inv, \Theta) \to \CC\UU_\ep(M_\ep, \B, \inv, \Theta)
\]
by $\tr_{\reg}$, the algebra $\UU_\ep(M_\ep, \B, \inv, \Theta)$ is very rarely free over $\CC\UU_\ep(M_\ep, \B, \inv, \Theta)$, so this is 
not a regular trace in the setting described in Sect. \ref{4.1}.
\ere
%%%%%%%%%%%%%
\sectionnew{Intersections of Cayley--Hamilton algebras and maximal orders}
%%%%
In this section, we start with background material on maximal orders and reduced traces, and then prove
a general theorem on intersections of maximal orders and Caylay--Hamilton algebras.
\subsection{The regular, standard and reduced traces} 
\label{traces}
\hfill

(1) If $R$ is a $\kk$-algebra over a commutative ring $\kk$, which is free over its center of finite rank $r$, 
then the left action of $R$ on itself induces the ring homomorphism 
\[
R \to \End_{\ZZ(R)}(R) \cong M_r(\ZZ(R)). 
\]
The composition of this map with the matrix trace $\Tr$ on $M_r(\ZZ(R))$ gives the {\em{regular trace}} of $R$:
\[
\tr_{\reg}: R \to \End_{\ZZ(R)}(R) \cong M_r(\ZZ(R)) \stackrel{\Tr}{\to} \ZZ(R).
\]

(2) If $R$ is a prime PI algebra of PI degree $d$ and $Q$ is the quotient field of $\ZZ(R)$, then we have the embedding 
\[
R \hra S:= R[(\ZZ(R)\backslash \{0\})^{-1}] \cong
R\otimes_{\ZZ(R)}Q,
\]
and by Posner's theorem \cite[Theorem 13.6.5]{MR}, $S$ is a central simple algebra of dimension $d^2$ with center $Q$. There is a finite field extension $F$ of $Q$ 
such that $S \otimes_Q F \cong M_d(F)$ (called a {\em{splitting field}} of $S$, see e.g. \cite[Sect. 7b]{Re}); that is $R\otimes_{\ZZ(R)}F \cong M_d(F)$.

The {\em{standard trace}} $\tr_{\sta} : R \to \ZZ(R) $ is the composition
\[
\tr_{\sta} : R \to R\otimes_{\ZZ(R)} F \cong M_d(F) \stackrel{\tr_{\reg}}{\lra} F. 
\]

The {\em{reduced trace}} $\tr_{\red} : R \to \ZZ(R) $ is the composition
\[
\tr_{\red} : R \to R\otimes_{\ZZ(R)} F \cong M_d(F) \stackrel{\Tr}{\lra} F. 
\]
Since the pair $(M_d(F), F, \Tr)$ is a Cayley--Hamilton algebra of degree $d$ if $\charr F \notin [1,d]$, the pair $(R, \ZZ(R), \tr_{\red})$ 
is a Cayley--Hamilton algebra if $\charr \ZZ(R) \notin [1,d]$. 

The standard and reduced traces are related by
\[
\tr_{\sta} = d \tr_{\red},
\]
cf. \cite[Eq. (2.5)]{BY}. 

(3) If we are in both situations (1) and (2), i.e., $R$ is a prime affine PI algebra which is free over its center of rank $r$, then the PI degree of $R$ is 
\begin{equation}
\label{PI-degree}
d = \sqrt{r}
\end{equation}
and 
\begin{equation}
\label{trace-ident}
\tr_{\reg} = \tr_{\sta} = d \tr_{\red}.
\end{equation}
The first equality follows at once from the fact that in this situation a $\ZZ(R)$-basis of $R$ gives an $F$-basis of $R \otimes_{\ZZ(R)} F$. 

When we need to emphasize the algebra $R$ in question, the above traces will be denoted by $\tr_{\reg}^R$, $\tr_{\sta}^R$ and $\tr_{\red}^R$.
%%%%%%%%%%
\def\fM{\mathfrak{M}}
\subsection{Intersections of maximal orders and Cayley--Hamilton algebras}
\bth{inter-max-order} 
Assume that 
\[
\{ R_\gamma \mid \gamma \in \Gamma \}
\]
is a collection of $\kk$-algebras over a commutative ring $\kk$,
which are maximal orders in a central simple algebra $S$ with center $Q$ of dimension $d^2$ (over $Q$) for an index set $\Gamma$.
For each $\gamma \in \Gamma$, let $R'_\gamma$ be a $\kk$-subalgebra of $R_\gamma$ and $Z_\gamma \subseteq \ZZ(R'_\gamma)$ 
be a central $\kk$-subalgebra of $R'_\gamma$ such that 
\[
\tr_{\red}^{R_\gamma}(R'_\gamma) \subseteq Z_\gamma.
\]
Set 
\[
R:= \bigcap_{\gamma \in \Gamma} R'_\ga
\quad \mbox{and} \quad
Z:= \bigcap_{\gamma \in \Gamma} Z_\ga \subseteq \ZZ(R).
\]
The following hold:
\begin{enumerate}
\item $\tr_{\red}^{R_\beta}|_R = \tr_{\red}^{R_\gamma}|_R$ for all $\beta, \gamma \in \Gamma$. Set $\tr_{\red} := \tr_{\red}^{R_\gamma}|_R$, 
which is independent on the choice of $\gamma \in \Gamma$.
\item $\tr_{\red}(R) \subseteq Z$.
\item If $\charr Q \notin [1,d]$, then the triple $(R, Z, \tr_{\red})$ 
is a Cayley--Hamilton algebra of degree $d$.
\end{enumerate}
For the rest of the theorem, we restrict ourselves to the case when $R'_\gamma = R_\gamma$, $Z_\gamma = \ZZ(R_\gamma)$ for all $\gamma \in \Gamma$.
\begin{enumerate}
\item[(4)] $Z$ is integrally closed. 
\item[(5)] If $S$ is a central localization of $R$, then $Z= \ZZ(R)$.
\item[(6)] If $R_\gamma$ is a central localization of $R$ then $R$ is a maximal order in $S$ whose reduced trace equals $\tr_{\red}$. 
\end{enumerate}
\eth
If $\charr Q \notin [1,d]$, then each triple $(R_\gamma, Z_\gamma, \tr_{\red}^{R_\gamma})$ is a Cayley--Hamilton algebra of degree $d$
by \leref{restrict-CH}. 
The third part of the theorem shows that their intersection is also a Cayley--Hamilton algebra of degree $d$. The sixth part of the 
theorem proves that an intersection of maximal orders is also a maximal order under the natural assumption that each
$R_\gamma$ is a central localization of $R$. 
\begin{proof}
(1) For all $\beta, \gamma \in \Gamma$, we have the commutative diagram
\[
\begin{tikzpicture}[scale=1.4]
\node (A) at (0,1) {$R$}; 
\node (B) at (1,1) {$R_\beta$};
\node (C) at (0,0) {$R_\gamma$};
\node (D) at (1,0) {$S$};
\draw
(A) edge[right hook- >,font=\scriptsize,] (B)
(A) edge[right hook ->,font=\scriptsize,>=angle 90] (C)
(B) edge[right hook ->,font=\scriptsize,] (D)
(C) edge[right hook ->,font=\scriptsize] (D);
\end{tikzpicture}
\]
By the construction of the reduced trace, the maps $\tr_{\red}^{R_\beta} : R_\beta \to \ZZ(R_\beta)$ and $\tr_{\red}^{R_\gamma} : R_\gamma \to \ZZ(R_\gamma)$ are the restrictions 
of $\tr_{\red}^S : S \to \ZZ(S)$ to the algebras $R_\beta$ and $R_\gamma$. Therefore, 
\[
\tr_{\red}^{R_\beta}|_R = \tr_{\red}^{S}|_R = \tr_{\red}^{R_\gamma}|_R.
\]

(2) Part (1) implies that $\tr_{\red}(R) \subseteq Z_\gamma$ for all $\gamma \in \Gamma$. Hence,  
\[
\tr_{\red}(R) \subseteq \bigcap_{\gamma \in \Gamma} Z_\gamma=Z.
\]

Part (3) follows from \leref{restrict-CH} and part (2). 

(4) Since $R_\gamma$ is a maximal order for all $\gamma \in \Gamma$,
$\ZZ(R_\gamma)$ is integrally closed (in $Q$). Therefore $Z = \bigcap_{\gamma \in \Gamma} \ZZ(R_\gamma)$ is 
integrally closed in $Q$. Thus, $Z$ is integrally closed in its quotient field, which is canonically identified with 
a subfield of $Q$.

(5) The assumption in part (5) implies that for all $\gamma \in \Gamma$, 
\begin{align*}
\ZZ(R) & = \{ z \in R \mid z r = rz, \forall r \in R\} \\
&=\{ z \in R \mid z s = s z, \forall s \in S\}  \\
&\subseteq \{ z \in R_\gamma \mid z s = s z, \forall s \in S\} = \ZZ(R_\gamma).  
\end{align*}
Hence, $\ZZ(R) \subseteq Z$, and by part (2), $\ZZ(R) = Z$.

(6) Let us prove that $R$ is an order. As $R_\gamma$ is a central localization of $R$ we have  $\ZZ(R) = \ZZ(R_\gamma) \cap R$. It follows that $\ZZ(R)$ is a domain.

Let us prove that  $\ZZ(R_\gamma) \subseteq \Fr(\ZZ(R))$, the field of fractions of $\ZZ(R)$. Let $z \in \ZZ(R_\gamma)$. By localization, we have $z= r u^{-1}$ where $r\in R$ and $0\neq u\in \ZZ(R)$. Since $z$ commutes with each element $R_\gamma$, it commutes with each element of $R$. It follows easily that $r$ commutes with each element of $R$. Hence $r\in \ZZ(R)$, and $z\in \Fr(Z(R))$. 

Consequently $Q= \Fr(\ZZ(R_\gamma) =\Fr(\ZZ(R))$.

Since $R\,  \Fr(\ZZ(R))$ contains  both $R_\gamma$ and $Q$, we have  $R\,  \Fr(\ZZ(R))\supset  R_\gamma\, Q = S$. Thus $R\,  \Fr(\ZZ(R))= S$. 
By \prref{r.maxorder}, the ring $R$ is an order of $S$.

Let us prove $R$ is a maximal order. Assume $0\neq z \in \ZZ(R)$ and a subring $\wt{R}$ of $S$ satisfy
\begin{equation}
R \subseteq \wt{R} \subseteq \frac1z  R.
\label{eq.1a}
\end{equation}
We will show that $\wt{R}=R$. Then by \prref{r.maxorder}, $R$ is a maximal order.

Fix $\gamma \in \Gamma$. There is a multiplicative set $\fM$ of central elements of $R$ such that $R \fM^{-1}= R_\gamma$. 
As $\wt{R} \subseteq \frac1z R$ and $z$ and each element of $\fM$ are central in $R$, each element of $\fM$ is central in $\wt{R}$ as well. 
In particular, one can define the localization $\wt{R} \fM^{-1}$. By localizing each algebra in \eqref{eq.1a} using the multiplicative set $\fM$, we get
$$ R_\gamma \subseteq \wt{R} \fM^{-1} \subseteq   \frac 1z  R_\gamma.$$
Since $R_\gamma$ is a maximal order, we conclude $\wt{R} \fM^{-1}=R_\gamma$. 
Because the elements of $\fM$ are invertible in $R_\gamma$, and thus are not zero divisors in $\wt{R}$, 
we have that $\wt{R} \subseteq \wt{R} \fM^{-1} = R_\gamma$. As this is true for all $\gamma \in \Gamma$, we have
$$\wt{R} \subseteq \bigcap_\gamma R_\gamma =R.$$
Therefore $\wt{R}= R$. Hence $R$ is a maximal order. Since each $R_\gamma$ is a localization of $R$, the reduced trace of $R$ equals $\tr_{\red}$. 
\end{proof}
%%%%%%%%%%%%%
\sectionnew{Root of unity upper quantum cluster algebras and maximal orders}
%%%%
In this section we show that all root of unity upper quantum cluster algebras are maximal orders, 
thus proving Theorem A from the introduction. We also prove part (5) of Theorem B from the introduction.
Throughout the section, we work over an arbitrary integral domain $\kk$ of characteristic 0.

\subsection{Trace maps on root of unity quantum tori}
\label{6.1}
Consider the root of unity quantum torus $\TT_{\ep}(\La)$ and denote
\[
\Ker(\La):=\{f\in \Zset^N \mid \La (f,g)=0\in \Zset/\ell, \forall g\in \Zset^N  \}. 
\]
$\Ker (\La)$ is a subgroup of $\Zset^N$ and thus a lattice. It has the same rank, i.e., a finite index in $\Zset^N$
\[
[\Zset^N :\Ker (\La)] < \infty
\]
because $\Ker (\La) \supseteq ( \ell \Zset)^N$. 

The center of $\TT_{\ep}(\La)$ is 
\begin{equation}
\label{cent-quant-tor}
\ZZ\big(\TT_{\ep}(\La)\big)=\Abbe\hspace{-0.1cm}-\hspace{-0.08cm}\Span \{ X^f \mid f \in \Ker(\La)\}.
\end{equation}

\bpr{quantum-torus-PI} The following hold for an arbitrary base ring $\kk$ which is an 
integral domain of characteristic 0:
\begin{enumerate}
\item $\TT_{\ep}(\La)$ is a free module over $\ZZ\big(\TT_{\ep}(\La)\big)$ of rank 
\[
[\Zset^N:\Ker(\La)]
\]
with basis $\{X^f \mid f \in \Delta\}$ where $\Delta \subset \Zset^N$ is a set of representatives for the cosets in $\Zset^N/\Ker(\La)$; 
\item $\TT_{\ep}(\La)$ and $\TT_{\ep}(\La)_{\geq}$ are maximal orders; 
\item If the base ring $\kk$ is an algebraically closed field containing the cyclotomic field $\Qset(\ep^{1/2})$, then $\TT_{\ep}(\La)$ is an Azumaya algebra.
\item The PI degree of $\TT_{\ep}(\La)$ equals
\[
d(\La):=\sqrt{[\Zset^N:\Ker(\La)]}.
\]
\end{enumerate}
\epr

\begin{proof}
Part (1) follows from the product formula \eqref{prod-quant-tor}. Part (2) follows from \cite[Theorem 6.5]{DP}. A more general statement will 
be proved in \thref{r.Maxorder} below. For part (3) see \cite[Proposition 7.2]{DP}.
Part (4) follows from \eqref{PI-degree} and part (1) of the proposition.
\end{proof}
\ble{reg-tr-quant-tor} The regular trace of $\TT_{\ep}(\La)$ is given by
\[
    \tr_{\reg}^{\TT_{\ep}(\La)}(X^f) = 
\begin{cases}
d(\La)^2\cdot X^f, & \text{if } f \in \Ker(\La)\\
0,                  &      \text{if } f\not\in \Ker(\La).
\end{cases}
\] 
\ele
The proof of the lemma is straightforward. 
\bco{red-tr-quant-tor} The reduced traces of $\TT_{\ep}(\La)$ and $\TT_{\ep}(\La)_\geq$ are given by
\[
\tr_{\red}^{\TT_{\ep}(\La)}(X^f) = 
\begin{cases}
d(\La)\cdot X^f, & \text{if } f \in \Ker(\La)\\
0,                  &      \text{if } f\not\in \Ker(\La).
\end{cases}
\]
\eco
The statement for $\TT_{\ep}(\La)$ follows from \eqref{trace-ident} and \leref{reg-tr-quant-tor}. The statement for $\TT_{\ep}(\La)_\geq$
follows from the fact that the reduced trace of $\TT_{\ep}(\La)_\geq$ is a restriction of that of $\TT_{\ep}(\La)$ because $\TT_{\ep}(\La)$ is a central localization 
of $\TT_{\ep}(\La)_\geq$.
\subsection{$\UU_\ep(M_\ep, \B, \inv, \Theta)$ is a maximal order}
Let us fix a subset $\Theta$ of the exchange graph $E_\ep(\La_{M_\ep},\B)$ of the root of unity upper quantum cluster algebra $\UU_\ep(M_\ep, \B, \inv)$. 
The center of the algebra $\UU_\ep(M_\ep, \B, \inv, \Theta)$ is given by 
\leref{centU}. The center $\ZZ(\TT_\ep(M'_\ep))$ of the quantum torus $\TT_\ep(M'_\ep)$ is explicitly described in \eqref{cent-quant-tor}. To determine 
the center of $\UU_\ep(M_\ep, \B, \inv, \Theta)$, we expand the elements of the algebra $\UU_\ep(M_\ep, \B, \inv, \Theta)$ in terms of the 
quantum torus $\TT_\ep(M'_\ep)$ for any seed $(M'_\ep, \B') \in \Theta$ (i.e., compute explicitly the embedding $\UU_\ep(M_\ep, \B, \inv, \Theta) \hra \TT_\ep(M'_\ep)$) 
and then apply \eqref{cent-quant-tor}. 

Since the quantum torus $\TT_\ep(M'_\ep)$ is a central localization of the algebra $\UU_\ep(M_\ep, \B, \inv, \Theta)$ (see Eq. \eqref{TU}) 
for any seed $(M'_\ep, \B')$ in $\Theta$, the PI degrees of the two algebras are equal to each other:

\bpr{PI-deg-U} The PI degree of $\UU_\ep(M_\ep, \B, \inv, \Theta)$ equals 
\[
d(\La_{M'_\ep}):=\sqrt{[\Zset^N:\Ker(\La_{M'_\ep})]}
\]
where $(M'_\ep, \B')$ is any seed in $\Theta$. 
\epr
As a consequence of the proposition, the integer $d(\La_{M'_\ep})$
is independent on the choice of a seed $(M'_\ep, \B')$ of $\UU_\ep(M_\ep, \B, \inv)$. We denote this integer by 
\[
d(M_\ep, \B).
\]

\bth{CH-U} Let $\kk$ be an integral domain of characteristic 0, $\ep^{1/2}$ be a primitive $\ell$-th root of unity 
in the algebraic closure of the faction field of $\kk$
(without any restriction on $\ell$) and $\UU_\ep(M_\ep, \B, \inv)$ be a root of unity upper quantum cluster algebra. 
Assume that $\Theta$ is a (not necessarily connected) subset of the exchange graph $E_\ep(\La_{M_\ep},\B)$ of $\UU_\ep(M_\ep, \B, \inv)$. The following hold:
\begin{enumerate}
\item For every pair of seeds $(M'_\ep, \B'), (M''_\ep, \B'') \in \Theta$, 
\[
\tr_{\red}^{\TT_\ep(M'_\ep)}\big{|}_{\UU_\ep(M_\ep, \B, \inv, \Theta)} = \tr_{\red}^{\TT_\ep(M''_\ep)}\big{|}_{\UU_\ep(M_\ep, \B, \inv, \Theta)}.
\]
We denote by $\tr_{\red}$ this restriction map coming from an arbitrary seed in $\Theta$. 
It satisfies $\tr_{\red} ( \UU_\ep(M_\ep, \B, \inv, \Theta)) \subseteq \ZZ( \UU_\ep(M_\ep, \B, \inv, \Theta))$.
\item The algebra $\UU_\ep(M_\ep, \B, \inv, \Theta)$ is a maximal order whose reduced trace is equal to $\tr_{\red}$. As a consequence, the triple 
\begin{equation}
\label{triple-aux}
(\UU_\ep(M_\ep, \B, \inv, \Theta), \ZZ(\UU_\ep(M_\ep, \B, \inv, \Theta)), \tr_{\red}) 
\end{equation}
is a Cayley--Hamilton algebra of degree equal to the PI degree $d(M_\ep,\B)$ of the algebra $\UU_\ep(M_\ep, \B, \inv, \Theta)$.
\item If the base ring $\kk$ is an algebraically closed field containing the cyclotomic field $\Qset(\ep^{1/2})$, then the union
\[
\bigcup_{(M'_\ep, \B') \in \Theta} \MaxSpec \big( \ZZ(\UU_\ep(M_\ep, \B, \inv, \Theta))[ M'_\ep(e_i)^{-\ell}, 1 \leq i \leq N] \big)
\]
inside $\MaxSpec \ZZ(\UU_\ep(M_\ep, \B, \inv, \Theta))$
is in the Azumaya locus of $\UU_\ep(M_\ep, \B, \inv, \Theta)$. 
\end{enumerate}
\eth
\begin{proof} 
Part (1) follows by applying parts (1-2) of \thref{inter-max-order} to the collection of mixed quantum tori
\[
\{ \TT_\ep(M'_\ep)_{\geq} \mid (M'_\ep, \B') \in \Theta \}
\]
which are maximal orders in the central simple algebra
\[
\TT_\ep(M_\ep)[(\ZZ(\TT_\ep(M_\ep)) \backslash \{0\})^{-1}]
\]
by \prref{quantum-torus-PI}(2). By definition, $\UU_\ep(M_\ep, \B, \inv, \Theta)$ is the intersection of these algebras. 
The triple in \eqref{triple-aux} is a Cayley--Hamilton algebra of the stated degree
by \exref{max-ord}. 

(2) By Eq. \eqref{TU} each $\TT_\ep(M'_\ep)_{\geq}$ is a central localization of $ \UU_\ep(M_\ep, \B, \inv, \Theta)$. 
Hence, by Theorem \ref{tinter-max-order}(6),  $ \UU_\ep(M_\ep, \B, \inv, \Theta)$ is a maximal order with reduced trace $\tr_{\red}$. 

Part (3) follows from the fact that each of the central localizations 
\[
\UU_\ep(M_\ep, \B, \inv, \Theta)[ M'_\ep(e_i)^{-\ell}, 1 \leq i \leq N]
\]
is isomorphic to a quantum torus (Eq. \eqref{TU}) which is an Azumaya algebra of PI degree $d(M_\ep,\B)$
(parts (2) and (4) of \prref{quantum-torus-PI} and \prref{PI-deg-U}).
\end{proof}
\bth{FullyAzumaya-U-CU} Assume the setting of part (3) of \thref{CH-U-CU}.
Then the union
\[
\bigcup_{(M'_\ep, \B') \in \Theta} \MaxSpec \big( \CC\UU_\ep(M_\ep, \B, \inv, \Theta)[ M'_\ep(e_i)^{-\ell}, 1 \leq i \leq N] \big)
\]
inside $\MaxSpec \CC \UU_\ep(M_\ep, \B, \inv, \Theta)$
is in the fully Azumaya locus of $\UU_\ep(M_\ep, \B, \inv, \Theta)$
with respect to its central subalgebra $\CC \UU_\ep(M_\ep, \B, \inv, \Theta)$,
cf. \deref{fullyAzumaya}. Over each such point, $\UU_\ep(M_\ep, \B, \inv, \Theta)$
has $\ell^N/ d(M_\ep, \B)$ irreducible representations of dimension $d(M_\ep, \B)$.
\eth
\begin{proof} Recall from \eqref{TU} that for every seed $(M'_\ep, \B') \in \Theta$, we have the isomorphism
\[
\UU_\ep(M_\ep, \B, \inv, \Theta)[ M'_\ep(e_i)^{-\ell}, 1 \leq i \leq N] \cong \TT_\ep(M'_\ep).
\]
Analogously one proves
\[
\CC\UU_\ep(M_\ep, \B, \inv, \Theta)[ M'_\ep(e_i)^{-\ell}, 1 \leq i \leq N] \cong \TT_\ep(M'_\ep)^\ell.
\]
The theorem follows from the two isomorphisms and the fact that 
the quantum torus in the right hand side of the first isomorphism is an Azumaya algebra of PI degree $d(M_\ep, \B)$
by \prref{quantum-torus-PI}(4).
\end{proof}
\bre{cluster-A-variety} Consider the case when $\Theta$ is the set of all seeds of the root of unity upper quantum cluster algebra 
$\UU_\ep(M_\ep, \B, \inv)$. The second isomorphism in \prref{cent-U} and the isomorphism of exchange graphs from Eq. \eqref{graph-iso} imply that the variety 
\[
\bigcup_{(M'_\ep, \B') \sim (M_\ep, \B)} \MaxSpec \big( \CC\UU_\ep(M_\ep, \B, \inv)[ M'_\ep(e_i)^{-\ell}, 1 \leq i \leq N] \big)
\]
in \thref{FullyAzumaya-U-CU} is isomorphic to the cluster $\mathcal{A}$-variety
\[
\bigcup_{(\wt{\mathbf{x}}', \B') \sim (\wt{\mathbf{x}}, \B)} \MaxSpec \kk[(x'_i)^{\pm 1}, 1 \leq i \leq N].
\]
\ere
\sectionnew{An example}

In the past it has been proved that quantum groups at roots of unity \cite{DK,DL,DP} as well as stated and Muller skein algebras at roots of unity of surfaces with nontrivial 
boundary  \cite{LY, Paprocki} are maximal orders in central simple algebras. These algebras are closely related to root of unity quantum cluster algebras with sufficiently 
large number of frozen variables and in all cases, the results were obtained by filtration arguments. On the other hand, 
root of unity quantum cluster algebras with fewer frozen variables behave in less tractable way and in general do not have appropriate filtrations. 
In this section we illustrate how our results apply to such an algebra.

\subsection{The root of unity quantum cluster algebra of type $A_2$} Consider the cluster algebra of type $A_2$ without frozen variables. 
It is associated to the exchange matrix
\[
\B 
= \begin{bmatrix} 
0 & 1
\\ -1 & 0
\end{bmatrix}.
\]
For a positive integer $\ell$ and a primitive $\ell$-th root of unity $\ep^{1/2}$, consider the root of unity quantum toric 
frame $M_\ep$ with matrix $\La = \ol{\B}$, where as in Sect. \ref{2.3}, the bar denotes reduction modulo $\ell$. 
The corresponding root of unity quantum cluster algebra and its upper counterpart will be denoted by 
$\AA_\ep(M_\ep, \B)$ and $\UU_\ep(M_\ep, \B)$, respectively; the inverted set $\inv$ is dropped from the 
notation since there are no frozen variables.

By \cite[Theorem 4.8]{NTY}, the exchange graph of $\AA_\ep(M_\ep, \B)$ is canonically isomorphic to the exchange graph of the underlying cluster algebra, which is well known to be a 
pentagon. The cluster variables of $\AA_\ep(M_\ep, \B)$ are
\begin{align*}
&X_1:= M_\ep(e_1), X_2:=M_\ep(e_2), Y_1:=M_\ep(-e_1 + e_2) + M_\ep(-e_1),  Y_2:=M_\ep(e_1 - e_2) + M_\ep(-e_2),
\\
&M_\ep(- e_1 - e_2) + M_\ep(-e_1) + M_\ep(-e_2).
\nn
\end{align*}
The proof of the following fact is direct and is left to the reader.
\bpr{A2} \hfill 
\begin{enumerate}
\item $\AA_\ep(M_\ep, \B) = \UU_\ep(M_\ep, \B)$ and $\AA_\ep(M_\ep, \B)$ equals the $\Abbe$-subalgebra of $\TT_\ep(M_\ep)$ 
generated by the first four cluster variables $X_1, X_2, Y_1, Y_2$. 
\item If $\ell$ is odd, then $\ZZ(\AA_\ep(M_\ep, \B) )= \CC\UU_\ep(M_\ep, \B)$.
\end{enumerate}
\epr
The first part of the proposition is in agreement with the results that for a cluster algebra with an acyclic and coprime seed, the cluster algebra equals its 
lower bound and the corresponding upper cluster algebra \cite[Theorem 1.18 and Corollary 1.19]{BFZ}. The last cluster variable equals 
\[
M_\ep(- e_1 - e_2) + M_\ep(-e_1) + M_\ep(-e_2) = \ep^{-1/2} Y_1 Y_2 - \ep^{-1/2}.
\]
\subsection{A presentation of the algebra} Consider the $\Abbe$-algebra $R$ with generators 
$x_1$, $x_2$, $y_1$, $y_2$ and relations
\begin{align*}
&x_2 x_1 = \ep^{-1} x_1 x_2,  &&y_2 x_1 = \ep x_1 y_2, &&&x_2 y_1 = \ep y_1 x_2,
\\
&x_1 y_1 = 1 + \ep^{1/2} x_2, &&y_1 x_1 = 1 + \ep^{-1/2} x_2, &&&
\\
&x_2 y_2 = 1 + \ep^{-1/2} x_1, &&y_2 x_2 = 1 + \ep^{1/2} x_1, &&&
\\
& y_2 y_1 = \ep^{-1} y_1 y_2 + (1 - \ep^{-1}).
\end{align*}
\bth{iso-max-ord} For all positive integers $\ell$ and primitive $\ell$-th roots of unity $\ep^{1/2}$
in the algebraic closure of the fraction field of the base ring $\kk$, which is assumed 
to be an integral domain of characteristic 0, the following hold:
\begin{enumerate} 
\item There is an $\Abbe$-algebra isomorphism $\theta : R \to \AA_\ep(M_\ep, \B)$ given by 
\begin{equation}
\label{theta}
\theta(x_k) = X_k, \; \; \theta(y_k) = Y_k \quad \mbox{for} \quad k=1,2.
\end{equation}
\item The algebra $R$ is a maximal order.
\end{enumerate}
\eth
\begin{proof}
(1) \prref{A2}(1) implies that $\AA_\ep(M_\ep, \B)$ is isomorphic to the subalgebra 
of the quantum torus
\[
\TT_\ep(M_\ep) = \frac{\Abbe \lcor X_1^{\pm 1}, X_2^{\pm 1} \rcor}{(X_2 X_1 - \ep^{-1} X_1 X_2)}
\]
generated by
\[
X_1, \; \; X_2, \; \; Y_1 = X_1^{-1} + \ep^{1/2} X_1^{-1} X_2, \; \; Y_2 = X_2^{-1} + \ep^{-1/2} X_2^{-1} X_1. 
\]
From this one easily verifies that the elements $X_1, X_2, Y_1, Y_2$ satisfy the defining relations for the generators of $R$. Thus, 
there is a well defined $\Abbe$-algebra homomorphism $\theta : R \to \AA_\ep(M_\ep, \B)$ given by \eqref{theta}. The defining relations 
of $R$ also imply that 
\[
R = \Abbe\hspace{-0.1cm}-\hspace{-0.08cm}\Span \{ y_1^{m_1} x_1^{n_1} x_2^{n_2} y_2^{m_2} 
\mid m_1, n_1, n_2, m_2 \in \Nset, \min(m_1, n_1) = \min (m_2, n_2) =0 \}.
\]
Using that $\{X_1^{n_1}, X_2^{n_2} \mid n_1, n_2 \in \Zset\}$ is an $\Abbe$-basis of $\TT_\ep(M_\ep)$, one easily shows that the set 
\[
\{ Y_1^{m_1} X_1^{n_1} X_2^{n_2} Y_2^{m_2} 
\mid m_1, n_1, n_2, m_2 \in \Nset, \min(m_1, n_1) = \min (m_2, n_2) =0 \}
\]
is an $\Abbe$-basis of $\AA_\ep(M_\ep, \B)$. Hence $\theta$ is an isomorphism.

Part (2) follows from the first part of the theorem and \thref{CH-U}(2). 
\end{proof}

\sectionnew{Monomial subalgebras of quantum tori and their intersections over subsets of seeds}
In this section we show that any monomial subalgebra of a quantum torus is Cayley--Hamilton 
and give a necessary and sufficient condition for when it is a maximal order.
We generalize the results from the previous section to prove that the 
intersection of monomial subalgebras over a collection of seeds 
of a root of unity quantum cluster algebra is always a Cayley--Hamilton algebra.
Throughout the base ring $\kk$ is assumed to be an integral domain of characteristic 0.

\subsection{Monomial subalgebras of root of unity quantum tori}\label{sec.qtorus} 
For a root of unity quantum torus $\TT_{\ep}(\La)$ and a subset
\[
\Phi \subseteq \Zset^N,
\]
denote 
\[
\BA(\Phi) := \Abbe\hspace{-0.1cm}-\hspace{-0.08cm}\Span \{ X^f \mid f \in \Phi \}. 
\]
Eq. \eqref{prod-quant-tor} implies that $\BA(\Phi)$ is an $\Abbe$-subalgebra of $\TT_{\ep}(\La)$ if an only if $\Phi$ is a 
submonoid of $\Zset^N$. 
\bde{monoidal-subl} A subalgebra of $\TT_{\ep}(\La)$ of the form $\BA(\Phi)$ will be called a {\em{monoidal subalgebra}}.
\ede
\bex{ex-monoid-suaalg}
\begin{enumerate}
\item Each mixed quantum torus $\TT_{\ep}(\La)_{\geq}$ is a monoidal subalgebra of $\TT_{\ep}(\La)$,
namely,
\[
\TT_{\ep}(\La)_{\geq} = \BA(\Phi) \quad \mbox{for} \quad
\Phi := \{(f_1, \ldots, f_N) \in \Zset^N \mid f_j \geq 0 ; \; \forall j \notin \ex \sqcup \inv\}.
\]
\item If $\Phi$ is a subgroup of $\Zset^N$, then $\BA(\Phi)$ is isomorphic to a root of unity (based) quantum torus for 
the restriction of the bicharacter $\La$ to $\Phi$.
\end{enumerate}
\eex
\bpr{monoid-CH} For every monomial subalgebra $\BA(\Phi)$ of a root of unity quantum torus $\TT_{\ep}(\La)$, the triple
\[
(\BA(\Phi), \BA(\Phi) \cap \ZZ(\TT_{\ep}(\La)), \tr_{\red}^{\TT_{\ep}(\La)})
\]
is a Cayley--Hamilton algebra of degree
\[
d(\La)=\sqrt{[\Zset^N:\Ker(\La)]},
\]
where $\tr_{\red}^{\TT_{\ep}(\La)}$ is the reduced trace of $\TT_{\ep}(\La)$ as in Sect. \ref{6.1}.
\epr
Since, $\ZZ(\BA(\Phi)) \subseteq \BA(\Phi) \cap \ZZ(\TT_{\ep}(\La))$, \prref{monoid-CH} implies that 
$(\BA(\Phi), \ZZ(\BA(\Phi)), \tr_{\red}^{\TT_{\ep}(\La)})$ is a Cayley--Hamilton algebra of the same degree. 
\begin{proof} It follows from \coref{red-tr-quant-tor} that 
\[
\tr_{\red}^{\TT_{\ep}(\La)}( \BA(\Phi)) \subseteq  \BA(\Phi) \cap \ZZ(\TT_{\ep}(\La)).
\]
The proposition now follows from \leref{restrict-CH} and parts (2) and (4) of \prref{quantum-torus-PI}.
\end{proof}

For every submonoid $\Phi \subseteq \Zset^N$, 
\[
\uL:= \Phi - \Phi= \{\al - \be \mid \al, \be \in \Phi\}
\]
is a subgroup of $\Zset^N$ and $\BA(\uL)$ is 
isomorphic to a root of unity quantum torus as described in \exref{ex-monoid-suaalg}(2). Furthermore, $\BA(\uL)$
is a central localization of $\BA(\Phi)$ because of  \eqref{prod-quant-tor} and the fact that each element of the lattice 
$\uL$ is of the form 
\[
\al - \be = \big( \al + (\ell-1) \be \big) - \ell \be
\]
for some $\al, \be \in \Phi$. Hence, the restriction of the reduced trace of the root of unity quantum torus $\BA(\uL)$ 
to $\BA(\Phi)$ equals the reduced trace of $\BA(\uL)$ (to be denoted by $\tr_{\red}^{\BA(\Phi)}$)
and the PI degrees of $\BA(\uL)$ and $\BA(\Phi)$ equal 
\begin{equation}
\label{PI-deg-A-Phi}
\sqrt{[\uL:\Ker(\La|_{\uL})]}.
\end{equation}
Applying \prref{monoid-CH} to the algebra $\BA(\Phi)$ and the quantum torus $\BA(\uL)$ leads to the following: 
\bpr{monoid-CH-2} For every monomial subalgebra $\BA(\Phi)$ of a root of unity quantum torus $\TT_{\ep}(\La)$, the triple
\[
(\BA(\Phi), \ZZ(\BA(\Phi)), \tr_{\red}^{\BA(\Phi)})
\]
is a Cayley--Hamilton algebra of degree \eqref{PI-deg-A-Phi}.
\epr
%%%%
\subsection{Criterion for a monomial subalgebra to be a maximal order} 
\label{monoid-max-ord}
Suppose $\Phi \subseteq \Zset^N$ is a submonoid, generating the group $\uL \subseteq \Zset^N$.  We consider $\Zset^N$  as a subset of $\Rset^N$.  

We say $\Phi$ is {\em integrally convex} if  $f,g\in \Phi$ implies that any point of $\uL$ lying in the interval segment $[f,g]$ is in $\Phi$. This is equivalent to the seemly weaker condition, 
which says that if $f\in \uL$ is such that $kf\in \Phi$ for a positive integer $k$ then $f\in \Phi$. In fact, assume the weaker condition and let $h\in \uL$ be in the segment $[f,g]$, 
where $f,g\in \Phi$. Then $h$ a non-negative rational linear combination of $f$ and $g$. Consequently there is a positive integer $k$ such that $kh$ is a non-negative integer linear combination of $f$ and $g$, and hence belongs to $\Phi$. Thus the weaker condition implies $h\in \Phi$.

For example the submonoid $\Nset \setminus \{1\}$ of $\Zset$ is not integrally convex.

For a non-zero $f \in\Rset^N$ the  {\em ray passing through $f$} 
is   $\Ray(f) =\{ t f \mid t\in \Rset , t >0 \}$. For a subset $V \subseteq \Rset^N$ let $\Ray(V)$ be the set of all rays passing through non-zero elements of $V$. There is a bijection $p: \Ray(\Rset^N) \to S$, where $S$ is the unit sphere in $\Rset^N$, given by
$ p(L) = L \cap S$.  We topologize $\Ray(\Rset^N)$ using the bijection $p$ and the topology of the unit sphere $S$.
We say $\Phi$ is {\em integrally closed} if $\Ray(\Phi)$ is closed in $\Ray(\uL)$. 

% For a ray $L$ let $L_\Zset$ be the set $L \cap \Zset^N$ of all integral points on $L$,  and if $L_\Zset \neq \emptyset$ we call $L$   {\em an integral} ray.

%Let $\Ray(\Phi)$ be the set of all rays passing through points in $\Phi$. We say $\Phi$ is {\em complete} if for any ray $L\in \Ray(\Phi)$ we have $L \cap \Phi= L \cap \uL$. In other words $\nexists$ $f \in \uL \backslash \Phi$ and $n \in \Zset_+$ such that $n f \in \Phi$.

%\no{We say $\Phi$ is {\em integrally closed} if an integral ray which is in the limit of $\Ray(\Phi)$ is in $\Ray(\Phi)$. Let us expand on the definition. The map $L\to p(L) = L \cap S$, where $S$ is the unit sphere in $\Rset^N$,  is a bijection between the set of all rays and $S$.   A point in $S$ is called {\em ray-integral} if it is the image of an integral ray. Let $C(\Phi)= p(\Ray(\Phi))$ and $\ol{C}(\Phi)$ be its topological closure. Then $\ol{C}(\Phi)$ is a convex subset of the unit sphere $S$, and $\Phi$ is integrally closed if and only if every ray-integral point in $\ol{C}(\Phi)$ is in $C(\Phi)$. }

We give below typical integrally closed and non-integrally closed submonoids.

\bex{cl-not-cl}
Suppose  $L_1, \dots, L_k$ are linear forms on $\Rset^N$, not necessarily having integral coefficients.
\begin{enumerate}
\item The submonoid determined by the non-strict inequalities
\[
\CL(L_1,\dots,L_k) := \{ f \in \Zset^N \mid L_i(f ) \ge 0, \ i=1,\dots, k\}
\]
is integrally convex and  integrally closed.
\item Consider the submonoid defined by the strict and non-strict inequalities 
$$\CL'(L_1,\dots,L_k) =  \{ f \in \Zset^N \mid L_1(f )> 0 \ \mbox{for} \ f \neq 0, \ L_i(f ) \ge 0, i=2, \dots, k \}.$$
Suppose $L_1$ has integer coefficients and is non-redundant in the definition of \\  $\CL'(L_1,\dots,L_k)$. Then $\CL'(L_1,\dots,L_k)$  is  integrally convex but not integrally closed.
\end{enumerate}
\eex

Note that integrally convex and integrally closed properties are intrinsic, meaning they do not depend on how the monoid $\Phi$ embeds in an abelian group $\Zset^N$.
\bth{r.Maxorder}
A monomial subalgebra $\BA(\Phi)$ of a root of unity quantum torus $\TT_{\ep}(\La)$ is a maximal order if and only if $\Phi$ is integrally convex and integrally closed.
\eth
For example, a very special case of the theorem shows that the mixed quantum tori $\TT_\ep(M_\ep)_\geq$ are maximal orders, 
see Examples \ref{eex-monoid-suaalg}(1) and \ref{ecl-not-cl}(1).

\begin{proof} By replacing $\Zset^N$ with $\uL$, we can assume that $\uL= \Zset^N$. Then the ring of fractions $\Fr(\BA(\Phi))$ of $\BA(\Phi)$ is equal to that of $\TT_{\ep}(\La)$. 
%\redc{The closed cone $\ol{C}(\Phi)$ is a subset of maximum dimension of the unit sphere, and has infinitely many  ray-integral interior points if $N >1$.}

(a) Assume that $\BA(\Phi)$ is a maximal order. 

Let us show that $\Phi$ is integrally convex. Assume $ f\in \Zset^N$ such that $kf\in \Phi$ where $k$ is a positive integer. We need to show that $f \in \Phi$. By replacing $k$ with a multiple of it we can assume $k$ is divisible by the order of $\ep$. Then $X^{kg}$ is central in $\BA(\Phi)$ for any $g\in \Phi$.
 Let $B$ be the algebra generated by $\BA(\Phi)$ and the monomial $b=X^f$. We will write $x =_\times y$ if $x=uy$ where $u$ is an invertible scalar. Since 
for any monomial $X^g$ we have $b X^g =_\times  X^g b$ and $\BA(\Phi)$ has a basis consisting of monomials, we have  $b \BA(\Phi) = \BA(\Phi) b$. It follows that
$$
B = \sum_{i=0}^\infty b^i \BA(\Phi)= \sum_{i=0}^{k-1} b^i \BA(\Phi),
$$
where the second identity follows from the fact that $b^k\in \BA(\Phi)$. Since $\Phi-\Phi= \Zset^N$, there are $f', f''\in \Phi$ such that $f'-f''=f$. Let $a= X^{f'}$ and $c= X^{f''}.$ Then for $ i=0, \dots, k-1$,
$$ b^i =_\times c^{-i} a^i  = c^{-k}  c^{k-i} a^i    \in c^{-k} \BA(\Phi).$$
It follows that $B \subseteq c^{-k} \BA(\Phi)$, and maximal order property implies $B= \BA(\Phi)$. Hence  $b \in B = \BA(\Phi)$, or $f\in \Phi$.  Thus $\Phi$ is integrally convex.

Let us now show that $\Phi$ is  integrally closed. If $N=1$ then the  integral convexity implies that $\Phi=\Nset$ or $\Phi= -\Nset$ or $\Phi=\Zset$. In each case $\Phi$ is integrally closed. 

Suppose now $N>1$.  Assume   $L\in \Ray(\Zset^N)$ is a limit of rays in $ \Ray(\Phi)$. We have to show that $L \in \Ray(\Phi)$. Let $C(\Phi) = p(\Ray(\Phi))$ and $C(\Zset^N) = p(\Ray(\Zset^N))$. The topological closure
 $\ol{C}(\Phi)$ of $C(\Phi)$  is a convex subset of $S$ of dimension equal to that of $S$.

Claim: If $x$ is an interior point of $\ol{C}(\Phi)$ and $x\in C(\Zset^N)$, then $x\in C(\Phi)$. This is because  $\Ray(x)$ is a convex linear combination of the rays 
in $\Ray(\Phi)$ when $x$ is an interior point of $\ol{C}(\Phi)$. As all rays involved are integral, the coefficients of the linear combination can be chosen to be non-negative rational numbers.  This implies that  $\Ray(x)$ contains a point in $\Phi$, or $x\in C(\Phi)$. 

In particular, if $p(L)$ is an interior point of $\ol{C}(\Phi)$, then $L\in \Ray(\Phi)$.
Consider the remaining case when $p(L)$ is on the boundary of $\ol{C}(\Phi)$. The interior of $\ol{C}(\Phi)$ is non-empty as $N>1$. Choose $f\in \Phi$ such that  $p(\Ray( f))$ is an interior point of 
$\ol{C}(\Phi)$ and let $a= X^{\ord(\ep) f}$ where $\ord(\ep)$ is the order of $\ep$. Note  that $a$ is in the center $Z(\BA(\Phi))$ of $\BA(\Phi)$, which is absolutely integrally closed. Let $g\in L\cap \Zset^N$ and $b= X^{\ord(\ep) g}$. Then $b$ is in the field of fractions $\Fr(Z(\BA(\Phi)))$.

From the claim  and the integral convexity, we have $f +  k g \in \Phi$ for all positive integers $k$. This is because
 $p(\Ray (f+k g))$ is an interior point of $\ol{C}(\Phi)$ and at the same time an element of $C(\Zset^N)$. Thus
$ab^k\in Z(\BA(\Phi))$ for all positive integers $k$. The absolutely integrally closed property implies $b \in Z(\BA(\Phi)) \subseteq \BA(\Phi)$, which means  $L$ is in $\Ray(\Phi)$. 

(b) Assume that $\Phi$ is integrally convex and integrally closed. We will prove that $\BA(\Phi)$ is a maximal order.

Suppose that $0\neq z\in Z(\BA(\Phi))$ and $B$ is an $\Abbe$-algebra such that 
\[
\BA(\Phi) \subseteq B \subseteq \frac 1z \BA(\Phi).
\]
We need to show  $B= \BA(\Phi)$. Let $b\in B$ be a non-zero element.
As $\TT_{\ep}(\La)$ is a central localization of $\BA(\Phi)$ and at the same time a maximal order, we have $b\in \TT_{\ep}(\La)$. See the proof of Theorem \ref{tinter-max-order}(6). 
%\redc{Expand here?}

As $0 \neq b \in \TT_{\ep}(\La)$, there is a non-empty finite set $\supp(b) \subseteq \Zset^N$ such that
$$ b = \sum_{f \in \supp(b)} c_f X^f, \quad 0 \neq c_f \in \Abbe.$$
Let $\Newt(b)$, known as the Newton polytope,  be the convex hull in $\Rset^N$ of $\supp(b)$. Note that $b\in \BA(\Phi)$ if and only if $(\Newt(b) \cap \Zset^N) \subseteq \Phi$.

For two non-zero $b_1, b_2 \in \TT_{\ep}(\La)$ we have 
\begin{equation}
\Newt(b_1 b_2) = \Newt(b_1) + \Newt(b_2):=  \{ f_1 + f_2 \mid f_1 \in \Newt(b_1), f_2 \in \Newt(b_2) \}.
\label{eq.Newt}
\end{equation}
%where  $$ \Newt(b_1) + \Newt(b_2) = \{ f_1 + f_2 \mid f_1 \in \Newt(b_1), f_2 \in \Newt(b_2) \}.$$
In fact Identity \eqref{eq.Newt} for commutative rings of Laurent polynomials is known \cite[Proposition~19.4]{Gru}, and the easy proof there carries over to  quantum tori. 
%In particular, if $\bk\in \Newt(b)$ then $t \bk \in \Newt(b^t)$ for any positive integer $t$.

For all positive integers $k$, the fact $b^kz \in \BA(\Phi)$ means $(k\, \Newt(b) + \Newt(z)) \cap \Zset^N \subseteq \Phi$. 
The integrally closed property implies the ray passing through any point $f\in\Newt(b) \cap \Zset^N$ is in  $\Ray(\Phi)$, 
and the integral convexity further implies that $f\in \Phi$. Thus $\Newt(b) \cap \Zset^N \subseteq \Phi$, or $b \in \BA(\Phi)$. 
This completes the proof that $\BA(\Phi)$ is a maximal order.
\end{proof}

\bre{rem1} 
The sufficient condition of Theorem \ref{tr.Maxorder} is a generalization of a result of the second author and J. Paprocki \cite{Paprocki}, and the proof is a generalization of the proof given there.
\ere

\bre{rem2} 
The stated skein algebra of a surface with non-empty boundary \cite{LY}  has an $\Nset$-filtration whose associated graded algebra is a monomial algebra associated 
to an integrally convex and integrally closed submonoid, and consequently the stated skein algebra of a surface with non-empty boundary is a maximal order, see Theorem 7.5 therein.

For a marked surface with non-empty boundary and  with at least one marked point on each boundary component, Muller \cite{Muller} defined a quantum cluster algebra. 
The second author and J. Paprocki showed that the Muller quantum cluster algebra has an $\Nset$-filtration whose associated graded algebra is a monomial algebra associated 
to an integrally convex and integrally closed submonoid, and consequently it is a maximal order, see \cite[Section 10]{Paprocki}.
\ere
\subsection{A cluster theoretic intersection of monomial algebras}
\bth{CH-U2} Assume that $\ep^{1/2}$ is a primitive $\ell$-th root of unity 
in the algebraic closure of the fraction field of the base ring $\kk$, which is assumed 
to be an integral domain of characteristic 0
(no restrictions on $\ell$), and $\UU_\ep(M_\ep, \B, \inv)$ is a root of unity upper quantum cluster algebra.
Let $\Theta$ be a (not necessarily connected) subset of the exchange graph $E_\ep(\La_{M_\ep},\B)$ of $\UU_\ep(M_\ep, \B, \inv)$
and $\BA_\Sigma(\Phi_\Sigma)$ be a monomial subalgebra of $\TT_\ep(M'_\ep)$ for each seed $\Sigma = (M'_\ep, \B') \in \Theta$
for some submonoids $\Phi_\Sigma \subseteq \Zset^N$.
Denote 
\[
\BA := \bigcap_{\Sigma \in \Theta} \BA_\Sigma(\Phi_\Sigma). 
\]

The following hold:
\begin{enumerate}
\item For every pair of seeds $(M'_\ep, \B'), (M''_\ep, \B'') \in \Theta$, 
\[
\tr_{\red}^{\TT_\ep(M'_\ep)}\big{|}_{\BA} = \tr_{\red}^{\TT_\ep(M''_\ep)}\big{|}_{\BA}.
\]
We denote by $\tr_{\red}$ this restriction map coming from an arbitrary seed in $\Theta$. 
\item $\tr_{\red} ( \BA ) \subseteq \cap_{\Sigma \in \Theta} \ZZ(\BA_\Sigma(\Phi_\Sigma))$.
\item The triple 
\[
(\BA, \cap_{\Sigma \in \Theta} \ZZ(\BA_\Sigma(\Phi_\Sigma)), \tr_{\red}) 
\]
is a Cayley--Hamilton algebra of degree equal to $d(M_\ep,\B)$.
%\item If all submonoids $\Phi_\Sigma$ of $\Zset^N$ are integrally convex and integrally closed, then 
\end{enumerate}
\eth
The theorem follows from applying parts (1-3) of \thref{inter-max-order} to the collection of Cayley--Hamilton algebras 
$(\BA_\Sigma(\Phi_\Sigma), \ZZ(\BA_\Sigma(\Phi_\Sigma)), \tr_{\red}^{\TT_\ep(M'_\ep))})$
for the seeds $\Sigma=(M'_\ep, \B') \in \Theta$ 
constructed in \prref{monoid-CH-2}. 
The theorem generalizes \thref{CH-U} in that intersections of mixed quantum tori over a collection of seeds are replaced 
with arbitrary intersections of monomial algebras. However, in that more concrete situation, \thref{CH-U} establishes a stronger result 
that the former intersections are maximal orders rather than just Cayley--Hamilton algebras.
%%%%%%%%%%%%%%%%%%%%%%%%%%%%%%%%%%%%%%%%%%%%%%%%%%%%%%%%%%%%%%%%%%%%%%%%%%%

%%%%%%%%%%%%%%%%%%%%%%%%%%%%%%%%%%%%%%%%%%%%%%%%%%%%%%%%%%%%%%%%%%%%%%%%%%%

%%%%%%%%%%%%%%%%%%%%%%%%%%%%%%%%%%%%%%%%%%%%%%%%%%%%%%%%%%%%%%%%%%%%%%%%%%%
\end{document}